\documentclass[11pt]{amsart}
\usepackage{amsmath,amsthm, amscd, amssymb, amsfonts, mathrsfs}
\usepackage[all]{xy}
\usepackage{color}
\usepackage{enumitem}
\usepackage{mathdots}%Please turn this on again after you finish
\usepackage{hyperref}

\usepackage{hhline} 

\usepackage[active]{srcltx}%{vista inversa}

\newcommand{\PSL}{\mathbf{PSL}}
\newcommand{\Gb}{\boldsymbol{G}}
\newcommand{\PSp}{\mathbf{PSp}}

\newcommand{\PSU}{\mathbf{PSU}}

\newcommand{\SU}{\mathbf{SU}}

\newcommand{\nsbgp}{\vartriangleleft}

\newcommand{\Jf}{{\tt J}}

\newcommand{\tu}{\mathtt{u}}
\newcommand{\tv}{\mathtt{v}}

\newcommand{\diag}{\operatorname{diag}}

\newcommand{\supp}{\operatorname{supp}}

\newcommand{\Inn}{\operatorname{Inn}}

\newcommand{\clas}[2]{\mathcal C(#1,#2)}

\newcommand{\kc}{\mathbb F_q}

\newcommand\toba{\mathfrak B }
\newcommand\bpi{\boldsymbol{\pi}}

\newcommand{\trid}{\triangleright}

\newcommand{\Fc}{{\mathcal F}}

\newcommand{\W}{{\mathcal W}}
\newcommand{\Zc}{{\mathcal Z}}

\newcommand{\dw}{\dot{w}}
\newcommand{\Mi}{M^{\mathcal F}}

\newcommand{\kk}{\Bbbk}%{}

\newcommand{\ku}{\mathbb C}
\newcommand{\Gc}{{\mathcal G}}

\newcommand{\Z}{\mathbb Z}
\newcommand{\N}{{\mathbb N}}
\newcommand{\I}{{\mathbb I}}
\newcommand{\J}{{\mathbb J}}

\def\Hb{\mathbb H}

\newcommand{\G}{{\mathbb G}}
\newcommand{\B}{{\mathbb{B}}}
\newcommand{\T}{{\mathbb{T}}}
\newcommand{\U}{\mathbb{U}}
\newcommand\V{\mathbb V}
\newcommand{\Pa}{\mathbb{P}}
\newcommand{\Le}{\mathbb{L}}

\newcommand{\syp}{\mathsf{Sp}}

\newcommand{\M}{{\mathbb M}}
\newcommand{\Q}{{\mathsf Q}}
\newcommand{\F}{{\mathbb F}}
\newcommand{\C}{{\mathcal C}}

\newcommand{\urad}{{\Q}}

\newcommand{\GL}{\mathbf{GL}}
\newcommand{\GU}{\mathbf{GU}}

\newcommand{\SL}{\mathbf{SL}}
\newcommand{\Sp}{\mathbf{Sp}}

\newcommand{\Ort}{\mathbf{O}}
\newcommand{\Fr}{\operatorname{Fr}}

\newcommand{\Vc}{{\mathcal V}}

\newcommand{\Oc}{{\mathcal O}}
\newcommand{\oc}{{\mathcal O}}

\newcommand{\Aut}{\operatorname{Aut}}

\newcommand\ad{\operatorname{ad}}
\newcommand\Ad{\operatorname{Ad}}
\newcommand\Gsc{\G_{\operatorname{sc}}}

\numberwithin{equation}{section}
\theoremstyle{plain}

\newtheorem{lema}{Lemma}[section]
\newtheorem{theorem}[lema]{Theorem}
\newtheorem{cor}[lema]{Corollary}

\newtheorem{prop}[lema]{Proposition}

\newtheorem{question-app}{Question}

\theoremstyle{definition}
\newtheorem{definition}[lema]{Definition}

\theoremstyle{remark}
\newtheorem{obs}[lema]{Remark}
\newtheorem{rmk}[lema]{Remarks}
\newtheorem{step}{Case}

\newtheorem{stepa}{Step}
\newtheorem{paso}{Step}

\newcommand\id{\operatorname{id}}

\newcommand\st{\mathbb S_3}

\newcommand\ac{\mathbb A_4}

\newcommand\as{\mathbb A_6}

\newcommand\Sim{\mathbb S}

\newcommand\sei{\mathbb S_6}

\def\pf{\begin{proof}}
\def\epf{\end{proof}}

\theoremstyle{remark}

\newcounter{tabla}\stepcounter{tabla}
\renewcommand{\thetabla}{\Roman{tabla}}

\begin{document}

\renewcommand{\baselinestretch}{1.2}

\thispagestyle{empty}
%\vspace*{2in}
\title[Nichols algebras over unipotent classes]
{Finite-dimensional pointed Hopf algebras\newline over finite simple groups of Lie type II. \newline 
Unipotent classes in symplectic  groups}

\author[N. Andruskiewitsch, G. Carnovale, G. A. Garc\'ia]
{Nicol\'as Andruskiewitsch, Giovanna Carnovale and\newline Gast\'on Andr\'es Garc\'ia}

\thanks{
This work was partially supported by
ANPCyT-Foncyt, CONICET, Secyt (UNC), the grant CPDA125818/12 of the University of Padova,
GNSAGA, the Visiting Scientist Program of the University of Padova
and the bilateral agreement between the Universities of C\'ordoba and Padova.}

\address{\hspace{-20pt}  N. A. : FaMAF, 
Universidad Nacional de C\'ordoba. CIEM -- CONICET. %\newline \noindent
Medina Allende s/n (5000) Ciudad Universitaria, C\'ordoba,
Argentina}
\email{andrus@famaf.unc.edu.ar}
\address{\hspace{-20pt} G. C.:
Dipartimento di Matematica,
Universit\`a degli Studi di Padova,
via Trieste 63, 35121 Padova, Italia}
\email{carnoval@math.unipd.it}
\address{\hspace{-20pt} G. A. G.: Departamento de Matem\'atica, Facultad de Ciencias Exactas,
Universidad Nacional de La Plata. CONICET. C. C. 172, (1900)
La Plata, Argentina.}
\email{ggarcia@mate.unlp.edu.ar}

%\subjclass[2010]{16T05}
%\date{\today}

\begin{abstract}
We show that Nichols algebras of most simple Yetter-Drin\-feld modules over the projective 
symplectic linear group over a finite field, 
corresponding to unipotent orbits, have infinite dimension. We give a criterium to deal with 
unipotent classes of general finite simple
groups of Lie type and apply it to regular classes.
\end{abstract}

\maketitle

%\begin{quote}{\textit{That is not dead which can eternal lie.\\
%And with strange aeons even death may die.} }\end{quote}
%\rightline{Abdul Alhazred}

\rightline{\href{https://www.youtube.com/watch?feature=player_embedded&v=sM0XvXEt2lo}
{\textit{The call of cthulhu}}}

\setcounter{tocdepth}{1}

\tableofcontents

\section{Introduction}

This is the second paper of a series intended to determine the finite-dimensional
pointed Hopf algebras with group of group-likes isomorphic to a finite simple group of Lie type.
An Introduction to the whole series was given in Part I \cite{ACG-I}. The base field is $\ku$.
Let $p$ be a prime number, $m\in \N$, $q =p^m$ and $\kc$ the field with $q$ elements.
In this paper  we consider Nichols algebras associated to unipotent conjugacy classes in
\emph{symplectic groups} $\Gb = \PSp_{2n}(q)$, $n \ge 2$, see e. g. \cite{W, MT}.
We consider also here the non-simple group $\PSp_4(2)\simeq\Sim_6$ for convenience.

Let  $\oc$ a conjugacy class of $\Gb$.
We seek to determine all $\Oc$ that \emph{collapse} \cite[2.2]{AFGV-ampa},
that is, the dimension of the Nichols algebra $\toba(\Oc, q)$ is infinite for
every finite faithful 2-cocycle $q$. Our main result says

\begin{theorem}\label{th:unipotent-chevalley-collapse}
Let $\Oc$ be a unipotent conjugacy class in $\Gb$.
If $\Oc$ is not listed in \ref{tab:uno}, then it collapses.
\end{theorem}

\begin{gather}\notag
\PSp_{2n}(q)  \\ \label{tab:uno} \tag*{Table \thetabla}
\begin{tabular}{|c|c|c|c|}
\hline $n$  & type  & $q$ & Remark  \\
\hline
$\geq 2$  & $W(1)^a\oplus V(2)$  & even   & cthulhu,  Lemma \ref{lem:sp-hook}
\\\cline{2-4}
&  $(1^{r_1}, 2)$  & odd, 9 or  &  cthulhu, Lemma  \ref{lem:sp-hook} \\
&    &  not a square &  
\\
\hline
$3$ & $W(1)\oplus W(2)$ & 2 & cthulhu, Lemma \ref{lem:W(2)W(1)not-D}\\
\hline
$2$ & $W(2)$  & even  & cthulhu,  Lemma  \ref{lem:sp2-(010)} \\
\cline{2-4}
& $(2, 2)$ & 3 &   one class cthulhu \\
&  &  & Lemma  \ref{lem:sp43}  \\
\cline{2-4}
& $V(2)^2$ & 2 & cthulhu, Lemma \ref{lem:22q2-1}  \\\hline
\end{tabular}
\end{gather}

In the next paper of the series we will deal with the non-semisimple classes in $\Gb$;
for this we will also need to consider the unipotent classes in the finite unitary groups.

\subsection*{Notation} We denote the cardinal of a set $X$ by $\vert X\vert$. 
If $k < \ell$ are positive integers, then we set $\I_{k, \ell} = \{i\in \N: k\le i\leq \ell\}$ and simply
$\I_{\ell} = \I_{1, \ell}$.

Let $G$ be a group; $N < G$, respectively $N \nsbgp G$, means that $N$ is a subgroup, 
respectively a normal subgroup, of $G$.
The centralizer, respectively the normalizer, of $x\in G$ is denoted by $C_{G}(x)$, 
respectively $N_{G}(x)$;
the inner automorphism defined by conjugation by $x$ is denoted by $\Ad x$. 
If $\Fc\in \Aut G$, then $G^{\Fc}$
denotes the subgroup consisting of points fixed by $\Fc$.

\section{Preliminaries on racks}\label{sec:racks-preliminaries}

Recall that a rack is a non-empty set $X$ with a self-distributive binary 
operation $\trid$ such that $x\trid \underline{\,\,\,}$ 
is bijective for all $x\in X$. The archetype of a rack is a conjugacy class in
a group with the conjugation operation. 
This notion allows considerable flexibility in the treatement of the conjugacy classes.

\subsection{Collapsing criteria}\label{subsec:racks-criteria}
We use criteria from \cite{ACG-I, AFGV-ampa} to prove 
Theorem \ref{th:unipotent-chevalley-collapse}; see \cite{ACG-I} for more details.
Let $X$ be a rack. One says that  

\medbreak
\noindent$\circ$ $X$ is \emph{of type D}
provided that there is a decomposable subrack
$Y = R\coprod S$  with elements $r\in R$, $s\in S$ such that
\begin{equation}\label{eq:typeD-rack}
r\trid(s\trid(r\trid s)) \neq s.
\end{equation}

\medbreak
\noindent $\circ$ $X$ is \emph{of type F} if it has a family of
mutually disjoint subracks $(R_a)_{a \in \I_4}$ and a family $(r_a)_{a \in \I_4}$ such that 
for all $a,b \in \I_4$
\begin{itemize}
  \item $R_a \triangleright R_b = R_b$;
  \item$r_a\in R_a$  and $r_a\triangleright r_b \neq r_b$ when $a\neq b$.
\end{itemize}

\medbreak
\noindent $\circ$ $X$ is \emph{cthulhu} if it is neither of type D nor of type F.

\medbreak
\noindent $\circ$ $X$ is \emph{sober} if every subrack is either abelian or indecomposable.

\begin{theorem}\label{th:type4} 
\cite[Theorem  3.6]{AFGV-ampa}; \cite[Theorem  2.2]{ACG-I}.
A rack $X$ of type D (respectively, F) collapses. \qed
\end{theorem}

\begin{lema}\label{lem:parabolic}Let $G$ be a finite group, 
$P < G$, $\pi\colon P\to L$ a quotient map and $x\in P$.
If $\Oc_{\pi(x)}^{L}$ is of type D, respectively F, then $\Oc_x^{G}$ is again so.
\end{lema}
\pf The class $\Oc_x^{G}$ contains the subrack $\Oc_x^{P}$ and $\pi$ 
induces a rack epimorphism $\Oc_x^{P}\to \Oc_{\pi(x)}^{L}$. 
The statement follows from \cite[Remark 2.9]{ACG-I}.\epf

Recall the convention in \cite{ACG-I}: all racks considered in this series are crossed sets.

\begin{lema}\label{lema:ACG-2-10}\cite[Lemma 2.10 (i)]{ACG-I}
Let $X$ and $Y$ be racks. Assume that there are $y_1\neq y_2\in Y$, $x_1\neq x_2\in X$
	such that  $x_1 \trid (x_2\trid (x_1\trid x_2)) \neq x_2$, $y_1\trid y_2 = y_2$.
Then $X \times Y$ is of type D. \qed
\end{lema}

\subsection{Conjugacy classes and subgroups}
\label{subsec:racks-remarks}
For recursive reasoning, we need to consider how a conjugacy 
class splits when intersected with a subgroup.
Let $G$ be a finite group,  $N < G$ and $x\in N$. Let
$\clas{N}{x}$ be the set of $N$-conjugacy classes contained in $\Oc_x^{G}$. 
We start with the case $N$ normal.

\begin{obs}\label{rem:conjug-classes-normal-subgps} 
\cite[Remark 2.1]{ACG-I}
If $N \nsbgp G$, then $\Oc_x^G$ is a union of $N$-conjugacy classes
isomorphic to each other as racks.
\end{obs}

Next relevant case is $N = G^{\Fc}$, where $\Fc\in \Aut G$.
Recall that $G$ acts on itself
by $x\rightharpoonup y = x y \Fc(x^{-1})$, $x, y \in G$.
Let $H^1(\Fc,G)$ be the set of $\Fc$-twisted conjugacy classes in $G$, i. e.
the orbits with respect to the action $\rightharpoonup$.

\begin{obs}\label{rem:general-gps} Let $M < G$ be $\Fc$-stable, 
$g,h\in G$, $x\in G^{\Fc}$. Set $z:= g^{-1}\Fc(g)$,  $w:= h^{-1}\Fc(h)$.
Then

\begin{enumerate}\renewcommand{\theenumi}{\alph{enumi}}
\renewcommand{\labelenumi}{(\theenumi)}
  \item\label{item:a} $gxg^{-1} \in  G^{\Fc} \iff  z\in C_G(x)$. 
  
  \item\label{item:b} $gMg^{-1}$ is $\Fc$-stable $\iff  z\in N_G(M)$.
  
  \item\label{item:c} Assume that $z\in N_G(M)$. Then $(gMg^{-1})^{\Fc} 
  = g(M^{\Ad z\circ\Fc})g^{-1}$.

  \item\label{item:d} Assume that $z\in C_G(M)$. Then $(gMg^{-1})^{\Fc} = g(M^{\Fc})g^{-1}$, 
  and hence $[(gMg^{-1})^{\Fc} ,(gMg^{-1})^{\Fc}] = g[\Mi,\Mi] g^{-1}$.
  
  \item\label{item:e} Assume that $z\in C_G(M)$ and $x\in M^{\Fc}$, respectively $x\in [\Mi,\Mi]$. 
  Then $gxg^{-1}\in (gMg^{-1})^{\Fc}$, respectively $\in [(gMg^{-1})^{\Fc}, (gMg^{-1})^{\Fc}]$, 
  and there are rack isomorphisms
\begin{align*}
\Oc_x^{\Mi} &\simeq \Oc_{gxg^{-1}}^{(gMg^{-1})^{\Fc}},& 
&\text{respectively}&   \Oc_x^{[\Mi ,\Mi]} 
&\simeq \Oc_{gxg^{-1}}^{[(gMg^{-1})^{\Fc}, (gMg^{-1})^{\Fc}]}.
\end{align*}

  \item\label{item:f} Assume that $z\in C_G(M)$, $x\in M^{\Fc}$ and $\Oc_x^{\Mi}$ 
  is of type D, respectively F,
  then $\Oc_{gxg^{-1}}^{G^{\Fc}}$ is of the same type.
  
  \item\label{item:g} Assume that $z\in C_G(M)$, $x\in[\Mi,\Mi]$ and $\Oc_x^{[\Mi,\Mi]}$ 
  is of type D, respectively F,
  then $\Oc_{gxg^{-1}}^{[G^{\Fc}, G^{\Fc}]}$ is of the same type. 

  \item\label{item:h} If  $z, w\in C_G(x)$, then $\Oc_{gxg^{-1}}^{G^{\Fc}} 
  = \Oc_{hxh^{-1}}^{G^{\Fc}}$ if and only if $z$ and $w$ belong 
  to the same $\Fc$-twisted conjugacy classes in $C_G(x)$.

\end{enumerate}
\end{obs}

\pf \eqref{item:a}, \eqref{item:b}, \eqref{item:c} and \eqref{item:d} 
are straightforward;  \eqref{item:e} 
follows from \eqref{item:a} and \eqref{item:d}, 
while  \eqref{item:f} and  \eqref{item:g}
follow from \eqref{item:e}.
\eqref{item:h}: Assume that $\Oc_{gxg^{-1}}^{G^{\Fc}} 
= \Oc_{hxh^{-1}}^{G^{\Fc}}$; take
$k\in G^{\Fc}$ such that $kgxg^{-1}k^{-1} = hxh^{-1}$ and 
$u = h^{-1}kg$. Then $u\in C_G(x)$ and $u\rightharpoonup z = w$.
Conversely, if $u\in C_G(x)$ satisfies $u\rightharpoonup z = w$, 
then $k = hug^{-1} \in G^{\Fc}$ and $kgxg^{-1}k^{-1} = hxh^{-1}$.
\epf

By Remark \ref{rem:general-gps} \eqref{item:h}, the map
$\varphi: \clas{G^{\Fc}}{x} \to H^1(\Fc, C_G(x))$ sending  
$\Oc_{gxg^{-1}}^{G^{\Fc}}$ 
to the class of $z$, is well-defined and injective.

\begin{lema}\label{lema:general-gps} Let $M < G$ be  $\Fc$-stable, 
such that $x\in M$. Assume that
every element in the image of $\varphi$ has a representative in $C_G(M) \subset C_G(x)$.

\smallbreak
\emph{(1)} If $x\in M^{\Fc}$  and $\Oc_x^{\Mi}$  is of type D, respectively F, 
then $\Oc$ is so for every $\Oc\in \clas{G^{\Fc}}{x}$.
 
\smallbreak
\emph{(2)} If $x\in [\Mi,\Mi]$  and $\Oc_x^{[\Mi,\Mi]}$  is of type D, respectively F,
then $\Oc$ is so for every $\Oc\in \clas{[G^{\Fc}, G^{\Fc}]}{x}$.
\end{lema}

\pf (1). By Remark \ref{rem:general-gps} \eqref{item:a} and the assumption,  
there exists $g\in G$ such that $gxg^{-1} \in  G^{\Fc}$,
$\Oc = \Oc_{gxg^{-1}}^{G^{\Fc}}$ and $z = g^{-1}\Fc(g)\in C_G(M)$.
Then  Remark \ref{rem:general-gps} \eqref{item:f} applies.
The proof of (2) is similar, using Remark \ref{rem:general-gps} \eqref{item:g}.
\epf

\section{Preliminaries on finite simple groups of  Lie type}
\label{sec:lie-preliminaries}

\subsection{Algebraic groups}\label{subsec:preliminaries-alg-gps} 
We mainly follow \cite{MT} as a source on algebraic groups 
and finite groups of Lie type,
with exceptions signaled along the text.

Let $\kk = \overline{\F_q}$ be the algebraic closure  of $\F_q$. 
All algebraic groups are affine and defined over $\kk$.
If $\Hb$ is an algebraic group, then $\Hb^\circ$ indicates the connected 
component of $\Hb$ containing the identity.
Also, $X(\Hb) = \text{Mor} (\Hb,\kk^{\times})$ is the group of characters  
of $\Hb$, and $X_{*}(\Hb) = \text{Mor} (\kk^{\times}, \Hb)$
is the set of multiplicative one-parameter subgroups in $\Hb$.

Let $\G$ be a simple algebraic group, $\G_{\ad}$ its adjoint quotient, $\Gsc$ 
its simply connected cover, with projection
$\bpi: \Gsc \to \G$.
We fix a maximal torus $\T$ of $\G$ and a Borel subgroup $\B$ containing it. 
The unipotent radical of $\B$ is denoted by $\U$.
We add a subscript ad or sc for the maximal torus and Borel 
of $\G_{\ad}$ or $\Gsc$;  our choices are compatible
with projections, e.~g. $\bpi(\T_{\operatorname{sc}}) = \T$.

The root system of $\G$ is denoted by $\Phi$, identified as
a subset of $X(\T)$;
the set of positive roots relative to $\T$ and $\B$ is denoted
by $\Phi^+$ and the simple roots by $\alpha_1,\,\ldots,\,\alpha_n$,
numbered as in \cite{Bou}.
The Weyl group $N_{\G}(\T)/\T$ is denoted by $W$; $(-,-)$ is
the $W$-invariant bilinear form on the ${\mathbb R}$-span of $\Phi$.
Let $\langle\,, \,\rangle: X(\T)\times X_{*}(\T) \to \mathbb{Z}$ be given by
$\langle\chi, \lambda\rangle = m$ if $(\chi \circ \lambda)(x) = x^{m}$.
The coroot system of $\G$ is denoted by 
$\Phi^{\vee} = \{\beta^{\vee}: \beta \in \Phi\} \subset X_{*}(\T)$,
where  $\langle\alpha, \beta^{\vee}\rangle = \frac{2(\alpha,\beta)}{(\beta,\beta)}$, 
for all $\alpha\in \Phi$.
Hence
\begin{align*}
\alpha(\beta^{\vee}(\zeta)) &= \zeta^{\frac{2(\alpha,\beta)}{(\beta,\beta)}}, 
& \alpha, \beta &\in \Phi, \zeta \in \kk^{\times}.
\end{align*}

For $\alpha\in\Phi$, there is a
monomorphism of abelian groups $x_\alpha: \kk\to \U$; we set
$\U_{\alpha}$ for the image of $x_\alpha$, called a root subgroup.
We adopt the normalization of $x_\alpha$ and the notation for 
the elements in $\T$ from \cite[8.1.4]{springer}.
We recall the commutation rule:
$t x_{\alpha}(a)t^{-1}=x_{\alpha}(\alpha(t)a)$, for $t\in \T$ and $\alpha\in\Phi$.
The group $\U$ is generated by the root subgroups
$\U_\alpha$, for $\alpha\in\Phi^+$.
More precisely, let us fix an arbitrary ordering on $\Phi^+$; then every
$u\in \U$ has a unique expression as a product (with respect to the fixed ordering)
\begin{align}\label{eq:unipotent-product}
u &= \prod_{\alpha\in\Phi^+}x_\alpha(c_\alpha),& c_\alpha &\in \kk, \alpha\in\Phi^+.
\end{align}
Let $\supp (u) = \{\alpha \in \Phi^{+}~|~
c_\alpha\neq 0\}$, that of course depends on the ordering.
In the sequel we will use frequently the Chevalley's commutator
formula \eqref{eq:chev} below, see \cite[Lemma 15, p. 22 and Corollary, p. 24]{yale}.
Let $\alpha, \beta \in\Phi^+$ such that $\alpha+\beta\in\Phi^+$. Fix a total order in the set
$\Gamma$ of pairs $(i, j)$ of positive integers such that $i\alpha+j\beta\in \Phi$.
Then there exist integers $c_{ij}^{\alpha\beta}$ such that
\begin{align}\label{eq:chev}x_{\alpha}(\xi)x_\beta(\eta)x_\alpha(\xi)^{-1}x_\beta(\eta)^{-1}
&= \prod_{(i,j)\in \Gamma}  x_{i\alpha + j\beta}(c_{ij}^{\alpha\beta}\xi^i\eta^j), &
\forall \xi, \eta&\in \kk.
\end{align}
(Clearly, \eqref{eq:chev} also holds when $\alpha+\beta$ is not a root, 
as $\U_\alpha$ and $\U_\beta$ commute in this case).
Let $m$, respectively $M$, be the maximum integer for which
$\beta-m\alpha\in\Phi$, respectively $\beta + M\alpha\in\Phi$.
Then the $\alpha$-string through $\beta$ is
the set of roots of the form $\beta-m\alpha$, \dots, $\beta + M\alpha$, 
and $m - M = \dfrac {2(\beta, \alpha)}{(\alpha, \alpha)}$.
It is known that, up to a nonzero scalar,  $c_{11}^{\alpha\beta}={m}+1$.
If the Dynkin diagram of $\G$ is simply-laced, then ${m}+1 = 1$;
otherwise, $\vert {m}+1\vert\in \{1, 2, 3\}$.
Then $c_{11}^{\alpha\beta}\neq0$ except in the cases listed in \ref{tab:c11-alpha-beta}.

\begin{gather} \label{tab:c11-alpha-beta} \stepcounter{tabla}\tag*{Table \thetabla}
\begin{tabular}{|c|c|c|c|}
\hline $p$  & type of $\Phi$ & $\alpha$ & $\beta$  \\
\hline  \hline
$3$  & $G_2$ & $\alpha_1$ & $2\alpha_1+\alpha_2$ \\
\cline{3-4}
& & $2\alpha_1 +\alpha_2$ & $\alpha_1$ \\
\cline{3-4}
& & $\alpha_1 +\alpha_2$ & $2\alpha_1+\alpha_2$ \\
\cline{3-4}
& & $2\alpha_1+\alpha_2$ & $\alpha_1 +\alpha_2$ \\
\hline \hline
$2$ & $B_n, C_n, F_4$  & \multicolumn{1}{r}{orthogonal} &to each other\\
\cline{2-4}
& $G_2$ & $\alpha_1$ & $\alpha_1+\alpha_2$ \\
\cline{3-4}
&  & $\alpha_1+\alpha_2$ & $\alpha_1$ \\
\hline
\end{tabular}
\end{gather}

\smallbreak
Let
 $\Sigma_{\alpha} = \{\beta \in \Phi^{+}: \alpha+\beta\in\Phi$ but  $(\alpha,\beta)$ 
 does not appear in  \ref{tab:c11-alpha-beta}$\}$,
 for $\alpha \in \Phi^{+}$. If $\beta \in \Sigma_{\alpha}$, then $x_{\alpha}(\xi)$ 
 and $x_{\beta}(\eta)$ do not commute  for  $\xi, \eta\in \kk^{\times}$.

\subsection{Conjugacy classes in finite simple groups of Lie type}
\label{subsec:preliminaries-gps}

\subsubsection{Finite simple groups of Lie type}\label{subsubsec:preliminaries-gps}
Let $\Hb$ be a semisimple algebraic group
defined over  $\kc$.  A
Steinberg endomorphism $F: \Hb\to\Hb$ is an abstract group
automorphism having a power equal to a Frobenius map
\cite[Definition 21.3]{MT}. We may assume that $F$ is the product of a Frobenius
endomorphism with an automorphism of $\Hb$
induced by a non-trivial Dynkin diagram automorphism.
The subgroup $\Hb^F$ is called a \emph{finite group of Lie type} \cite[Definition 21.6]{MT}.

\medbreak
Let $\G$ be a simple algebraic group and let $F$ 
be a Steinberg endomorphism of $\Gsc$.
Assume that it descends to a  Steinberg endomorphism 
of $\G$ (again called $F$), that happens 
when $\ker \bpi$ is $F$-stable,
see \cite[Example 22.8]{MT} for precise conditions. In particular, 
$F$ descends to $\G_{\ad} \simeq \G/Z(\G)$ 
always, and to  $\G$ when it is $\F_q$-split.
It is well-known that $\G_{\ad}$ is a simple abstract group \cite[Proposition 12.5]{MT}
but $\G_{\ad}^F$ is not simple in general.
However $\Gb := \Gsc^F/Z(\Gsc^F)$ is a  finite simple group
except for the following 8 examples  \cite[Theorem 24.17]{MT}:

\begin{itemize}
 \item $\PSL_2(2)\simeq \st$; $\PSL_2(3)\simeq \ac$; $\PSp_4(2)\simeq \sei$;
 \item  $\PSU_3(2)$, ${}^2B_2(2^{2})$ (both solvable);
 \item  $G_2(2) \simeq \Aut \PSU_3(3)$,
${}^2G_2(3) \simeq \Aut \PSU_2(8)$ (almost simple);
 \item ${}^2F_4(2)$, that contains a normal subgroup isomorphic
to the Tits group, with index 2.
\end{itemize}

Henceforth we assume that $\Gb = \Gsc^F/Z(\Gsc^F)$ 
is not one of these 8 groups and
call it a  \emph{finite simple group of Lie type}.
Notice that $\bpi (\Gsc^F) = [\G^F, \G^F]$ and there is
 the alternative useful  description
$\Gb \simeq [\G^F, \G^F]/ \bpi(Z(\Gsc^F))$. In particular,
$\Gb \simeq [\G_{\ad}^F, \G_{\ad}^F]$ \cite[Proposition 24.21]{MT}.

\subsubsection{Conjugacy classes}\label{subsubsec:preliminaries-classes}
There is a huge literature on the description of the conjugacy
classes in $\Gb$, see for instance the bibliography in  \cite{Hu, LS, MT}.
We shall give precise references as they are needed.
To start with, we recall the following arguments:

\smallbreak
$\diamond$ Every $F$-stable $\G$-conjugacy class  $\oc$
meets $\G^F$ \cite[Theorem 21.11 (a)]{MT}, a consequence of the
Lang-Steinberg Theorem \cite[Theorem 21.7]{MT}.

\smallbreak
$\diamond$ Let $\oc$ be an $F$-stable $\G$-conjugacy
class, $x\in \oc\cap \G^F$ and
\begin{equation}\label{eq:a(x)}
 A(x):= C_{\G}(x)/C_{\G}(x)^\circ.
\end{equation} 
Then $\clas{\G^F}{x}$ is
in bijection with  $H^1(F, A(x))$  \cite[8.5]{Hu}, \cite[Theorem 21.11 (b)]{MT}.

\smallbreak
From the preceding two facts, we see that to determine
the conjugacy classes in $\Gb$, one possible way is to 
consider the following questions:

\begin{enumerate} \renewcommand{\theenumi}
{\alph{enumi}}\renewcommand{\labelenumi}{(\theenumi)}
  \item\label{question:F-stable-conj-classes}
  Describe the $F$-stable $\G$-conjugacy classes.
  \item\label{question:GF-conj-classes} For a given $F$-stable
  $\G$-conjugacy class $\Oc$, describe the $\G^F$-conjugacy classes in $\oc\cap \G^F$.
  \item\label{question:pass-to-G} Pass this information to $\Gb$.
\end{enumerate}

\smallbreak These questions were treated in extent in the literature. 
We will recall the known answers for different kinds
of conjugacy classes along the way. Now we state some other useful facts.

\smallbreak
$\diamond$
The Borel subgroup $\B$ and the maximal torus $\T$ are
chosen $F$-stable, which is possible by \cite[Corollary 21.12]{MT}. Hence so is $\U =[\B, \B]$.

\smallbreak
$\diamond$ The subgroup $W^F$ of $F$-fixed points in 
the Weyl group $W$ is isomorphic to
$N_{\G^F}(\T)/ \T^F$ by \cite[Proposition 23.2]{MT}; clearly, 
$N_{\G^F}(\T) = N_{\G}(\T)\cap \G^F$. Even more,
every element in $W^F$ has a representative in 
$N_{[\G^F, \G^F]}(\T)$ \cite[Corollary 24.2]{MT}.

\subsection{Unipotent classes in  finite simple groups of Lie type}
\label{subsec:preliminaries-unipotent}

We need to describe the unipotent conjugacy classes
in finite simple groups of Lie type.
We keep the notations and assumptions from \ref{subsubsec:preliminaries-gps}
for $\G$, $F$ and $\Gb$;
let $\pi: \Gsc^F \to \Gb = \Gsc^F/Z(\Gsc^F)$ be the natural projection.
Every $x\in \Gsc$ has a Chevalley-Jordan decomposition
$x = x_sx_u= x_ux_s$, with $x_{s}$ semisimple and $x_{u}$ unipotent.
This decomposition boils down to the group $\G$ and to the finite groups $\G^F$, $[\G^F,\G^F]$
and $\Gb$, where it agrees with the decomposition
in the $p$-part, namely $x_u$, and the $p$-regular part, namely $x_s$.

\subsubsection{Unitary groups} \label{subsubsec:unitary}
In some inductive arguments we use the unitary groups $\PSU_n(q)$.
When dealing with them we will use the following matrix description.
Let ${\Jf}_n =\left(\begin{smallmatrix}
&&1\\
&\iddots&\\
1&&
\end{smallmatrix}\right) = {\Jf}_n^{-1} \in \GL_n(\kk)$. Let
$\Fr_q$, respectively $F$, be the Frobenius endomorphism of $\GL_n(\kk)$ 
raising all entries of the matrix to the $q$-th power,
respectively given by $F(X)={\Jf}_n\;^t\!(\Fr_q (X))^{-1}{\Jf}_n$, $X\in\GL_n(\kk)$. 
Following \cite[Examples 21.14(2), 23.10(2)]{MT}, the unitary and special unitary groups are
$\GU_n(q) = \GL_n(\kk)^F$, $\SU_n(q) = \SL_n(\kk)^F$.
Also, $\SU_n(q)$ can be realized as a subgroup of $\SL_n(q^2)$ \cite{W}. 
If $h\in \N$, then
\begin{align*}
F^{2h}(X) &= \Fr_{q^{2h}}(X),& F^{2h+1}(X) &= 
{\Jf}_n\;^t\!(\Fr_{q^{2h + 1}} (X))^{-1}{\Jf}_n.
\end{align*}
Hence $\GU_n(q)$, respectively $\SU_n(q)$, $\PSU_n(q)$,
can be identified with a subgroup of $\GU_n(q^{2h+1})$, 
respectively $\SU_n(q^{2h+1})$, $\PSU_n(q^{2h+1})$.

\smallbreak
The unipotent conjugacy classes in $\SU_n(q)$ are described as the
unipotent conjugacy classes in $\SL_n(q)$.
Indeed,

\smallbreak
\noindent $\diamond$
Every  unipotent class in $\SU_n(q)$ has a type : $u \in \SU_n(q)$
is of type $\lambda = (\lambda_{1},\ldots, \lambda_{k})$
if the elementary
factors of its characteristic polynomial equal
$(X -1)^{\lambda_{1}}$, $(X -1)^{\lambda_{2}}$, \dots,
$(X -1)^{\lambda_{k}}$,
where $\lambda_{1}\geq \lambda_{2} \geq \cdots \geq \lambda_{k}$. Conversely, 
since all unipotent classes in $\G = \SL_n(\kk)$ are $F$-stable (by a direct computation), for any type 
 there is a unipotent class in $\SU_n(q)$, by the Lang-Steinberg theorem.

\smallbreak
\noindent  $\diamond$  By \cite[8.5]{Hu} every  unipotent class in $\GL_n(\kk)$ 
meets $\GU_n(q)$ in exactly one class,
 since $C_{\GL_n(\kk)}(x)$ is connected for every $x$ \cite[I.3.5]{sp-st}.

\smallbreak
\noindent  $\diamond$ Since $\SU_n(q)$ is normal in $\GU_n(q)$, 
Remark \ref{rem:conjug-classes-normal-subgps} says that 
all  unipotent classes in $\SU_n(q)$ with the same type are isomorphic as racks.

\subsubsection{The isogeny argument} 
Section \ref{sec:symplectic} is devoted to unipotent classes in Chevalley groups.
By the isogeny argument, Lemma \ref{lem:isogeny-argument} below,
it is enough to treat
the unipotent classes in $\Gsc^F$ or  $[\G^F,\G^F]$.
This takes care of Question \eqref{question:pass-to-G} in \ref{subsubsec:preliminaries-classes}
and gives  flexibility to choose $\G$ in a  suitable form, e.~g. in matrix form.
Let $\Gc$ be a semisimple algebraic, resp. finite, group and
$\Gc_u$ the set of unipotent, resp. $p$-elements, in $\Gc$.

\begin{lema}\label{lem:isogeny-argument} \cite[Lemma 1.2]{ACG-I}
Let $\Zc$ be a central subgroup of $\Gc$
whose elements are all semisimple, respectively $p$-regular.
Then the quotient map $\pi: \Gc \to \Gc/\Zc$
induces a rack isomorphism  $\pi: \Gc_u \to (\Gc/\Zc)_u$
and a bijection between the sets of $\Gc$-conjugacy classes
in $\Gc_u$ and  in $(\Gc/\Zc)_u$. \qed
\end{lema}

\subsubsection{A reduction argument}

The determination of the unipotent conjugacy classes in
$\G$ and those that are $F$-stable, 
Question \eqref{question:F-stable-conj-classes} in \ref{subsubsec:preliminaries-classes}, 
is well-known, see \cite[Chapter 7]{Hu}, \cite[Chapters 7, 17, 22]{LS}.  
But the description of the $\G^F$-conjugacy classes in $\oc\cap \G^F$ 
for an $F$-stable $\G$-conjugacy class $\Oc$, 
Question \eqref{question:GF-conj-classes} in \ref{subsubsec:preliminaries-classes},
is more delicate; for example there is a
class in $\Sp_4(\kk)$ which splits into 2 classes in $\Sp_4(q)$
of different size \cite[Table 8.1]{LS}. Similar examples occur for other groups.
This was not the case when $\G = \SL_n(\kk)$ and $F$ is $\F_q$-split by 
Remark \ref{rem:conjug-classes-normal-subgps}.
To start with, observe that  $\U^F$, which is
a $p$-Sylow subgroup of $\G^F$ and $[\G^F, \G^F]$ \cite[Corollary 24.11]{MT},
is isomorphic to its image in $\Gb$ by Lemma \ref{lem:isogeny-argument}.
Hence, every unipotent element in $\G^F$, or in $\Gb$, 
is conjugated to an element in $\U^F$.
Also,  the $[\G^F,\G^F]$-classes into which a $\G^F$-class
in $[\G^F,\G^F]$ splits are all isomorphic as racks,  
see Remark \ref{rem:conjug-classes-normal-subgps}.

\smallbreak
The following result combines Remark \ref{rem:general-gps},
Lemma \ref{lema:general-gps} and \cite[8.5]{Hu}.

\begin{lema}\label{lem:J}
Let $\M$ be an $F$-stable subgroup of $\G$ and $x\in \M^F$.

Let $g\in \G$ such that $z= g^{-1}F(g)\in C_{\G}(x)$.
Assume that the class of $z$ in $H^1(F, A(x))$ has a representative in $C_{\G}(\M)$.

\smallbreak
\begin{enumerate}\renewcommand{\theenumi}{\alph{enumi}}
\renewcommand{\labelenumi}{(\theenumi)}
  \item\label{item:J-a}  If $\Oc_x^{\M^{F}}$ is of type D, respectively F, 
  then $\Oc_{gxg^{-1}}^{\G^{F}}$ is so.
  
\smallbreak
  \item\label{item:J-b} If $x\in[\M^{F},\M^{F}]$ and $\Oc_x^{[\M^{F},\M^{F}]}$ 
  is of type D, respectively F,
  then $\Oc_{gxg^{-1}}^{[\G^{F}, \G^{F}]}$ is so. 
   
\end{enumerate}

\smallbreak
   Assume that every element in $H^1(F, A(x))$ has a representative in $C_{\G}(\M)$.
This happens for instance if $C_{\G}(x) = C_{\G}(\M) \Hb$ with $\Hb$ connected.
  
\smallbreak
 \begin{enumerate}\renewcommand{\theenumi}{\alph{enumi}}
 \renewcommand{\labelenumi}{(\theenumi)} \setcounter{enumi}{2}
  \item\label{item:J-c}If  $\Oc_x^{\M^{F}}$  is of type D, resp. F,
  then $\Oc$ is so for every $\Oc\in \clas{\G^{F}}{x}$.
    
\smallbreak
  \item\label{item:J-d} If $x\in [\M^{F},\M^{F}]$  and $\Oc_x^{[\M^{F},\M^{F}]}$  
  is of type D, respectively F, then $\Oc$ is so for every $\Oc\in \clas{[\G^{F}, \G^{F}]}{x}$. \qed
  
\end{enumerate}

\end{lema}

\subsection{Criteria to collapse for unipotent classes} 
Let $\Gb$ be a finite simple group of Lie type and 
$\Oc$ a unipotent conjugacy class  in $\Gb$.
We realize $\Oc$ as a unipotent conjugacy class  in $[\G^{F},\G^{F}]$, 
where as above, $\G$ is a simple algebraic group and $F$ 
is a Steinberg endomorphism of $\G$.

\begin{definition}\label{def:alfabeta2}
Let $\alpha,\beta \in \Phi^+$ such that $\alpha+\beta\in\Phi$ 
but the pair $\alpha,\beta$ does not appear in   \ref{tab:c11-alpha-beta}.
We fix an ordering of $\Phi^+$.
We say that $\Oc$  has the
{\em $\alpha\beta$-property} if there exists $u\in \Oc\cap \U^F$ 
such that  $\alpha, \beta\in \supp u$ and
\begin{align}\label{eq:ab-property}
\begin{aligned}
 \alpha+\beta = \sum_{1\le i \le r} \gamma_i, &\text{ with } 
 r > 1, \gamma_i\in \supp u \\ &\implies r =2,\quad \{\gamma_1, \gamma_2\} = \{\alpha, \beta\}.
\end{aligned}
\end{align}
\end{definition}

\begin{rmk}\label{rem:phi-u} Let $u\in \Oc\cap \U^F$.

(i) If  there exist {\em simple} roots $\alpha$
and $\beta \in \supp u$ adjacent in the Dynkin diagram 
of $\Phi$ (so that $\alpha + \beta$ is a root),
then $\Oc$ has the $\alpha\beta$-property.

(ii) Let $\alpha, \beta \in \Phi^+$ such that $\Oc$ has the
$\alpha\beta$-property. By \eqref{eq:ab-property}, neither $\alpha$ nor $\beta$ can be
decomposed as a sum of roots in $\supp(u)$. Using the 
Chevalley commutator formula \eqref{eq:chev}, we infer that
$\alpha$ and $\beta$ lie in the support of $u$ for every ordering on
$\Phi^+$, and the $\alpha\beta$-property is independent of the  ordering.
\end{rmk}

\subsubsection{Unipotent classes of type D in Chevalley and Steinberg groups}
We give a criterium to determine if unipotent classes in
Chevalley and Steinberg groups are of type D.
\emph{In this subsection, we assume that $q$ is odd.} 
Recall that the only groups corresponding to 
very twisted Steinberg endomorphisms in odd characteristic 
are the Ree groups ${}^2G_2(3^{2h+1})$.
See \cite[Section 12.4]{carter-simple}.

\begin{prop}\label{prop:alfabetachevstein}
Let $\Gb$ be a finite simple group of Lie type. Assume
$\Oc$ has the $\alpha\beta$-property, for some $\alpha,\beta \in 
\Phi^{+}$ such that $q>3$
when $(\alpha,\beta)= 0$. Then  $\Oc$  is of type D.
\end{prop}

\pf

\begin{paso}\label{lem:alfabetaDsplit2}
If there exists
 $t\in \T\cap[\G^F,\G^F]$ such that $1\neq\alpha(t)\neq \beta(t)$, then $\Oc$ is of type D.
\end{paso}

Fix an ordering of the positive roots
ending with  $\alpha+\beta<\beta<\alpha$. Since $\Oc$
has the $\alpha\beta$-property, there exists
$u\in \Oc$ with
$$u=\prod_{\gamma\in\supp(u)}x_\gamma(a_\gamma) \in
\Bigg(\prod_{\substack{\gamma\in\supp(u)\\
\gamma\neq \alpha,\beta,\alpha+\beta}}\U_\gamma\Bigg)x_{\alpha+\beta}(a_{\alpha+\beta})
x_{\beta}(a_\beta)x_{\alpha}(a_\alpha),$$
and $a_{\alpha}a_{\beta}\neq 0$. Let $r=u$, $s=trt^{-1}\in \Oc$.
Then 
$$\langle r,\,s\rangle \subseteq H :=\langle \U_\gamma~|~\gamma\in\supp(u)\rangle.$$

Also $s\in\left(\prod_{\gamma\in\supp(u),\
\gamma\neq \alpha}\U_\gamma\right)x_{\alpha}(\alpha(t)a_\alpha)$; 
we see using \eqref{eq:chev} and \eqref{eq:ab-property} that
\begin{align*}
\Oc_r^H &\subseteq \Bigg(\prod_{\substack{
\delta =\gamma_{1}+\cdots+\gamma_{l}\\
\delta\neq\alpha,\ \gamma_i\in\supp(u)}}\U_\delta\Bigg)x_\alpha(a_\alpha), &
\Oc_s^H &\subseteq \Bigg(\prod_{\substack{
\delta =\gamma_{1}+\cdots+\gamma_{l}\\
\delta\neq\alpha,\ \gamma_i\in\supp(u)}}\U_\delta\Bigg)x_\alpha(\alpha(t)a_\alpha).
\end{align*}
Since $\alpha(t)\neq 1$, $\Oc_r^{H}\neq\Oc_s^{H}$, hence
$\Oc_r^{\langle r,\,s\rangle}\neq\Oc_s^{\langle r,\,s\rangle}$.
Since $rs,\,sr\in \U^F$ and $p\neq2$,
$(rs)^2\neq(sr)^2$ if and only if $rs\neq sr$. 
To prove the last inequality, and conclude that $\Oc$ is of type D,
let $V:=\langle \U_\gamma ~|~\gamma\in\supp(u),\, 
\gamma\neq \alpha,\beta,\alpha+\beta\rangle$.
Observe that if a right coclass $Vw$ of some $w\in \U$ 
contains an element of the form $x_{\alpha+\beta}(z)x_\beta(y)x_\alpha(x)$,
then $x,y,z$ are unique by \eqref{eq:ab-property} and
 Remark \ref{rem:phi-u} (ii), using \eqref{eq:chev}.
Again by \eqref{eq:ab-property} and Remark \ref{rem:phi-u} (ii), 
using \eqref{eq:chev}, we see that the coclass $Vrs$ contains
$x_{\alpha+\beta}(z)x_\beta((1+\beta(t))a_\beta)x_\alpha((1+\alpha(t))a_\alpha)$, with
\begin{align*}
z = \beta(t) c_{11}^{\alpha,\beta}  
a_\alpha a_\beta+(1+(\alpha+\beta)(t))a_{\alpha+\beta},
\end{align*}
while $Vsr$ contains $x_{\alpha+\beta}(z')x_\beta((1+\beta(t))
a_\beta)x_\alpha((1+\alpha(t))a_\alpha)$, with
\begin{align*}
z' = \alpha(t)c_{11}^{\alpha,\beta} a_\alpha a_\beta + 
(1+(\alpha+\beta)(t)) a_{\alpha+\beta}.
\end{align*}
Since
$\alpha(t)\neq\beta(t)$ and $c_{11}^{\alpha,\beta}a_\alpha a_\beta  \neq 0$ 
by assumption, we get that $rs\neq sr$.

\begin{paso}\label{paso:alfabetaD-chevalley}
If $\Gb$ is a Chevalley group, then there exists
 $t\in \T\cap[\G^F,\G^F]$ such that $1\neq\alpha(t)\neq \beta(t)$.
\end{paso}

Without loss of generality, if $\alpha$ and $\beta$ have different lengths,
we choose $\beta$ to be the longest one.
Take $t=\beta^\vee(\zeta)\in\T$, where $\zeta$ is a generator  of 
$\F_q^\times$. Then $t\in [\G^F,\,\G^F]$ by \cite[8.1.4]{springer},
and $ \alpha(t)= \zeta^{\frac{2(\alpha,\beta)}{(\beta,\beta)}}$, 
$\beta(t)=\zeta^2$. If $\Phi$ is simply-laced then 
$r= \frac{2(\alpha,\beta)}{(\beta,\beta)}=-1$
and $1 \neq \alpha(t)\neq  \beta(t)$. If $\Phi$ is of type $G_{2}$, then 
$r\in \{-1,1\}$ and the same assertion follows. If $\Phi$
is doubly-laced, then $r\in\{-1,0\}$. But if $r=0$, then 
$\beta(t) = \zeta^{2}\neq 1$ since by assumption $q>3$. Thus 
$1 \neq \beta(t)\neq \alpha(t) $ and
the claim follows by interchanging $\alpha$ and $\beta$. 

\smallbreak
Hence the Proposition for Chevalley groups follows from Steps \ref{lem:alfabetaDsplit2} 
and \ref{paso:alfabetaD-chevalley}.

\begin{paso}\label{paso:alfabetaD-Steinberg}
If $\Gb$ is a Steinberg group, then there exists
 $t\in \T\cap[\G^F,\G^F]$ such that $1\neq\alpha(t)\neq \beta(t)$.
\end{paso}
Here
$\Phi$ is simply-laced so $\frac{2(\alpha,\beta)}{(\beta,\beta)}=-1$. Assume first that the
Dynkin diagram automorphism $\theta$ associated with $F$ is an involution.
Then  the $\langle \theta\rangle$-orbit of $\beta$ is either
$\{\beta\}$ or $\{\beta,\theta(\beta)\}$.
In the former case, take
$t=\beta^\vee(\zeta) \in\T$  for  a
generator $\zeta$ of $\F_q^\times$ and
conclude as in Step \ref{paso:alfabetaD-chevalley}.
In the latter,  take
$t=\beta^\vee(\xi)(\theta\beta)^\vee(\xi^q) \in \T$ for
a generator $\xi$ of $\F_{q^2}^\times$. 
Then 
$t\in \bpi(\G_{sc}^F)=[\G^F,\G^F]$; $\alpha(t)\in\{\xi^{-1},\xi^{-1\pm q},\,\xi^{-1+2q}\}$
and $\beta(t)\in\{\xi^2,\,\xi^{2-q}\}$. Hence $\alpha(t)\neq 1,\,\beta(t)$ unless $q=3$ and either
\begin{align}
\label{eq:caso2}& (\alpha,\beta)=-1, &&(\alpha,\theta\beta)=0  &\text{ and }(\beta,\theta\beta)=-1, \text{ or }\\
\label{eq:caso1}& (\alpha,\beta)=-1, &&(\alpha,\theta\beta)=1  &\text{ and }(\beta,\theta\beta)=0.
\end{align}
Let $q=3$. If $\Phi$ is of type $D_n$ or $E_6$, then case \eqref{eq:caso2} never occurs because 
 $(\beta,\theta\beta)=0$ whenever $\beta\neq\theta\beta$. If $\Phi$ is of type $A_n$,
then case \eqref{eq:caso2} occurs only if  $\beta=\varepsilon_i-\varepsilon_j$, for $i<j$ and either: 
$\alpha=\varepsilon_l-\varepsilon_i$ for $l<i$, $2j=n+2$, and $l\neq j, n+2-i$,  or 
$\alpha=\varepsilon_j-\varepsilon_l$ for $j<l$,  $2i=n+2$; and $l\neq i, n+2-j$.
In both situations we take $t=\alpha^\vee(\xi)(\theta\alpha)^\vee(\xi^3)$ for $\xi$
a generator of $\F_9^\times$. This gives the claim in case \eqref{eq:caso2}.

By applying $\theta$ we observe, using Remark \ref{rem:phi-u} (ii)  that if $\alpha$ 
and $\beta$ satisfy condition \eqref{eq:ab-property}, then  $\theta\alpha$ and $\theta\beta$ also lie in $\supp(u)$. Therefore, \eqref{eq:ab-property} forces
$\theta(\alpha+\beta)\neq\alpha+\beta$.

If $\Phi$ is of type $A_n$, then a pair of roots satisfies \eqref{eq:caso1} only if $\beta=\varepsilon_i-\varepsilon_j$ with $\{i,j\}\cap\{n-i+2, n-j+2\}=\emptyset$ and either
$\alpha=\varepsilon_j-\varepsilon_{n-i+2}$ or $\alpha=\varepsilon_{n-j+2}-\varepsilon_i$ with $\left |\{i,j,n-j+2,n-i+2\}\right |=4$. 
Since in this case $\alpha+\beta$ would be $\theta$-invariant, such pairs are discarded.

If $\Phi$ is of type $D_n$, case \eqref{eq:caso1}
occurs only if $\beta=\varepsilon_i\pm\varepsilon_n$, 
$\alpha=\varepsilon_j\mp\varepsilon_n$ with $n\neq i\neq j\neq n$. We discard such pairs as we did for type $A_n$.
%Let $k\in \I_n\setminus\{i,j,n\}$ 
%The element $t=(\varepsilon_k-\varepsilon_j)^\vee(-1)$ will do.

Let $\Phi$ be of type $E_6$.  If a pair $(\alpha,\beta)$ satisfies \eqref{eq:caso1} 
and $(\beta,\alpha)$ does not, we interchange $\alpha$ and $\beta$. We verify by inspection that there are no 
pairs of roots $\alpha$ and $\beta$ such that  $(\alpha,\beta)$ and $(\beta,\alpha)$ are both in case \eqref{eq:caso1} and such that 
$\theta(\alpha+\beta)\neq\alpha+\beta$. This gives the claim when $\theta^2=1$.

\smallskip

Assume now $\Phi$ is of type $D_4$ and $\theta$ has order $3$.
We will show that $\alpha_2\in\{\alpha,\beta\}$. Let us fix an ordering of the roots in increasing height ${\rm ht}$ and let 
$u\in\mathcal{O}\cap \U^F$ be as in Definition \ref{def:alfabeta2}.  
We consider the support of $u$ with respect to this ordering.
The outer automorphism $\theta$ of order $3$ permutes $\alpha_1,\alpha_3$
and $\alpha_4$ and fixes $\alpha_2$. By inspection, for simple roots we have $\alpha\in \supp(u)$ 
if and only if $\theta(\alpha)\in\supp(u)$.  In addition, $\gamma+\gamma'\not\in\Phi$ if ${\rm ht}(\gamma)={\rm ht}(\gamma')\geq2$ 
or if ${\rm ht}(\gamma)={\rm ht}(\gamma')=1$ and $\gamma,\gamma'\neq\alpha_2$. 
So, if $\{\alpha_1,\alpha_2\}\not\subset\supp (u)$ we have $\alpha \in \supp(u)$ 
if and only if $\theta(\alpha) \in \supp (u)$ for every $\alpha\in\Phi^+$. 
Thus, if $\alpha_2\not\in\supp(u)$, condition (\ref{eq:ab-property}) is not verified 
for any pair $\alpha,\beta\in\supp(u)$ such that $\alpha+\beta\in \Phi$. 
So, $\alpha_2\in\supp(u)$. 
If $\alpha_1\in\supp(u)$ then we take $\alpha=\alpha_2$, $\beta=\alpha_1$. If, instead, $\alpha_1\not\in\supp(u)$, then $u$ 
has the $\alpha\beta$-property if and only if $\alpha_1+\alpha_2+\alpha_3+\alpha_4\in\supp(u)$ 
and $\{\alpha_1+\alpha_2,\alpha_2+\alpha_3+\alpha_4\}\not\subset\supp(u)$. 
In this case, we have $\alpha=\alpha_2$, $\beta=\alpha_1+\alpha_2+\alpha_3+\alpha_4$. 
In both cases, we take  $t=\alpha_1^\vee(\xi)(\theta\alpha_1)^\vee(\xi^q)(\theta^2\alpha_1)^\vee(\xi^{q^2})$
for $\xi$ a generator of $\F_{q^3}^\times$. 
Then, $t\in \bpi(\G_{sc}^F)=[\G^F,\G^F]$, and $\alpha(t)=\xi^{-(1+q+q^2)}$ and $\beta(t)\in\{\xi^2,\xi^{(1+q+q^2)}\}$. 
Then, $1\neq \alpha(t)\neq\beta(t)$ unless
$q=3$ and $\beta=\alpha_1+\alpha_2+\alpha_3+\alpha_4$. In this case, we replace $t$ by $t\alpha_2^\vee(-1)$.

%Then  the $\langle \theta\rangle$-orbit of $\beta$ is either
%$\{\beta\}$ or $\{\beta,\theta(\beta),\theta^2(\beta)\}$.
%The former case is treated as above. Analogously,
%we may assume that $\theta(\alpha) \neq \alpha$,
%since otherwise we interchage $\alpha$ and $\beta$.
%Thus for the latter case we assume that the
%$\langle \theta\rangle$-orbit of $\alpha$ and $\beta$ have both order $3$. Take
%$t=\beta^\vee(\xi)(\theta\beta)^\vee(\xi^q)(\theta^2\beta)^\vee(\xi^{q^2})$
 %for  a generator $\xi$ of $\F_{q^3}^\times$.  Then $t\in \bpi(\G_{sc}^F)=[\G^F,\G^F]$,
%$\alpha(t)=\xi^m$ where $m\in\{-1-q+q^2, -1+q-q^2,-1+q+q^2\}$
%and $\beta(t)=\xi^2$, since $(\beta,\theta\beta)=0$.
%As $0<|m|<q^3-1$, we have that 
%$\alpha(t)\neq1$ and $\xi^{m-2}\neq1$ so $\alpha(t)\neq\beta(t)$.

\smallbreak
Hence the Proposition for Steinberg groups follows from 
Steps \ref{lem:alfabetaDsplit2} and \ref{paso:alfabetaD-Steinberg}.

\begin{paso}\label{paso:alfabetaD-Ree}
If $\Gb = {}^2G_2(3^{2h+1})$, $h \ge 1$, then there exists
 $t\in \T\cap[\G^F,\G^F]$ such that $1\neq\alpha(t)\neq \beta(t)$. 
\end{paso}

In this case the possible (unordered)  pairs $\{\alpha,\beta\}$ are
\begin{align*}
&\{\alpha_1,\alpha_2\},& &\{\alpha_1,\alpha_1+\alpha_2\},& &\{\alpha_2,3\alpha_1+\alpha_2\}.
\end{align*}
The last two pairs are interchanged by the non-standard graph automorphism 
$\theta$ such that $x_{\alpha_1}(\zeta)\mapsto x_{\alpha_2}(\zeta^3)$ and 
$x_{\alpha_2}(\zeta)\mapsto x_{\alpha_1}(\zeta)$ for every $\zeta\in\kk$, \cite[12.4]{carter-simple}. 
Applying the Steinberg endomorphism $\Fr_{3^h}\circ\theta$ to a representative in a class $\Oc$
we see that $\Oc$ has the $\alpha\beta$-property for  $\{\alpha_1,\alpha_1+\alpha_2\}$ 
if and only if it has it for $\{\alpha_2,3\alpha_1+\alpha_2\}$. So, it is enough to consider  
$\{\alpha_1,\alpha_2\}$ and  $\{\alpha_1,\alpha_1+\alpha_2\}$.
Let $\zeta$ be a generator of $\F_{3^{2h+1}}^\times$ and let
$t=\alpha_1^\vee(\zeta^{3^h})\alpha_2^\vee(\zeta)$. Then $t\in \T^F\cap [\G^F,\G^F]$ and
%\begin{align*}
%1\neq\alpha_1(t)=\zeta^{2-3^h} &\neq \alpha_2(t)=\zeta^{2\cdot 3^h-3}, \\
%\alpha_1(t)=\zeta^{2-3^h} &\neq (\alpha_1+\alpha_2)(t)=\zeta^{3^h-1}.
%\end{align*}
\begin{align*}
1\neq\alpha_1(t)=\zeta^{2\cdot 3^h-1} &\neq \alpha_2(t)=\zeta^{-3^{h+1}+2}, \\
\alpha_1(t)= \zeta^{2\cdot 3^h-1} &\neq (\alpha_1+\alpha_2)(t)=\zeta^{3^h-1}.
\end{align*}
%\smallbreak
Hence the Proposition for Ree groups follows from 
Steps \ref{lem:alfabetaDsplit2} and \ref{paso:alfabetaD-Ree}.
\epf

\subsubsection{Unipotent classes of type F in Chevalley and Steinberg groups.}
In this subsection we address the case when 
$q$ is even albeit some results are valid more generally.
We give criteria to determine when a unipotent class is of type 
F in Chevalley or Steinberg groups.

\begin{prop}\label{prop:alfabetaF}
Assume that one of the following conditions hold:
\begin{itemize}
 \item $\Gb$ is a Chevalley group and $q\notin\{2,3,4,5, 7\}$;
 \item $\Gb=\PSU_3(q)$ and $q\not\in\{ 2,5,8\}$;
\item $\Gb$ is a Steinberg group and $q>8$.
\end{itemize}

If $\Oc$ is a unipotent class in $\Gb$ and has the $\alpha\beta$-property,
for some $\alpha,\beta \in \Phi^{+}$, then it is of type F.
\end{prop}

\pf 
\begin{stepa}\label{lem:alfabetaF}
If there exists a family $(t_a)_{a \in \I_4}$ in $\T\cap[\G^F,\G^F]$ such that
\begin{align}
\label{eq:alfabetaF2}
\alpha(t_a)\beta(t_b)&\neq\alpha(t_b)\beta(t_a)& &\mbox{ for every } a\neq b,
\end{align}
then $\Oc$ is of type F.
\end{stepa}

Notice that \eqref{eq:alfabetaF2} implies
\begin{align}
\label{eq:alfabetaF1}
(\alpha(t_a),\beta(t_a))&\neq(\alpha(t_b),\beta(t_b))& &\mbox{ for every }a\neq b.
\end{align}

Let $r_a:=t_a u t_a^{-1}$ and $R_a:=\U^F\triangleright r_a$,  $a\in \I_4$. 
We claim that \eqref{eq:alfabetaF2} ensures $r_a\triangleright r_b \neq r_b$ 
for every $a\neq b$,
and that \eqref{eq:alfabetaF1} ensures that $R=\coprod_{a\in \I_4}R_a$
 is a subrack with $R_a\triangleright R_b = R_b$.

As in Step ~\ref{lem:alfabetaDsplit2} of Proposition \ref{prop:alfabetachevstein}, we
fix an ordering of $\Phi^+$ ending with  $\alpha+\beta < \beta < \alpha$. 
%If $\alpha$ and $\beta$ have
%different lengths, we assume that $\alpha$ is the longest one.
Let $\V = \langle \U_\gamma ~|~\gamma\in\supp(r),\,
\gamma\neq \alpha,\beta,\alpha+\beta\rangle$.
Since $\Oc$ has the $\alpha\beta$-property, there exists
$r\in \Oc$ with
$r\in \V x_{\alpha+\beta}(a_{\alpha+\beta}) x_{\beta}(a_\beta)x_{\alpha}(a_\alpha)$
and $a_{\alpha}a_{\beta}\neq 0$.
Then
$r_{a}\in \V x_{\alpha+\beta}((\alpha+\beta)(t_a)a_{\alpha+\beta})
x_\beta(\beta(t_a)a_{\beta})x_\alpha(\alpha(t_a)a_{\alpha})
$.
By \eqref{eq:ab-property} and Remark \ref{rem:phi-u} (ii), 
using \eqref{eq:chev}, we see that the coclass $\V r_a r_b$ contains
$x_{\alpha+\beta}(x)
x_\beta(y)
x_\alpha(z)$
with 
\begin{align}
&x=(\alpha+\beta)(t_a) a_{\alpha+\beta} + c_{11}^{\alpha,\beta}
\alpha(t_a)\beta(t_b)a_{\beta}a_{\alpha},\\
&y=(\beta(t_a) + \beta(t_b))a_{\beta},\\
&z=(\alpha(t_a)+\alpha(t_b))a_{\alpha}.
\end{align} 
Arguing as in the proof of Step ~\ref{lem:alfabetaDsplit2} of 
Proposition \ref{prop:alfabetachevstein}, we see that
$r_a\trid r_b\neq r_b$ and that
$R_a\subset \V' x_\beta(\beta(t_a)y)x_\alpha(\alpha(t_a)x)$ with
$\V' = \langle \U_\gamma~|~\gamma\neq \alpha,\beta\rangle$. 
By a direct computation; $\coprod_{a\in \I_4} R_a$ is subrack of 
$\Oc$, hence $\Oc$ is of type F.

\begin{stepa}\label{lem:alfabetaF-chev}
 If $\Gb$ is a Chevalley group and $q>7$, then there exists a family $(t_a)_{a \in \I_4}$ in $\T\cap[\G^F,\G^F]$ 
satisfying \eqref{eq:alfabetaF2}.
\end{stepa}
Let $\zeta$ be a generator of $\F_q^\times$. By assumption on $q$, 
for $e_a:=a-1$ with $a\in \I_4$ we have
$re_a\not\equiv re_b\,\mod\; (q-1)$ for all pairs $a\neq b$ and $1\leq r \leq 3$. 
If $\alpha$ and $\beta$ have different lengths, we assume that $\alpha$ is the longest one.
Set $t_a=\alpha^\vee(\zeta^{e_a}) \in\T$; by \cite[8.1.4]{springer}, 
$t_a\in [\G^F,\G^F]$. 
Then $\alpha(t_a)\beta(t_b)=\zeta^{2e_a+me_b}$
with $m \in \{-1, 0 ,1\}$ and a 
direct verification gives the claim.

\medbreak
Hence the Proposition for Chevalley groups follows from 
Steps \ref{lem:alfabetaF} and \ref{lem:alfabetaF-chev}.

\begin{stepa}\label{lem:alfabetaF-psu} 
If $\Gb= \PSU_3(q)$, for $q\not\in \{2,5,8\}$, then there exists 
$t_a\in \T\cap[\G^F,\G^F]$ for $a\in \I_4$ satisfying \eqref{eq:alfabetaF2}.
\end{stepa}

The only classes with the $\alpha\beta$-property are the regular ones. 
In this case we have $\alpha=\alpha_1$, $\beta=\theta(\alpha_1)=\alpha_2$ 
and, for $\zeta$ a generator in $\F_{q^2}^\times$ we set
\begin{equation}\label{eq:ta}t_a=\alpha^\vee(\zeta^{a-1})
\beta^\vee(\zeta^{(a-1) q}),\quad\mbox{ for }a\in \I_4.\end{equation}

Then  $(\alpha(t_a),\beta(t_a))=(\zeta^{(a-1)(2-q)},\zeta^{(a-1)(2q-1)})$
and the claim follows from a direct computation.

\medbreak
Hence the Proposition for $\PSU_3(q)$ follows from Steps \ref{lem:alfabetaF} 
and \ref{lem:alfabetaF-psu}.
We assume in the remaining Steps   \ref{lem:alfabetaF-ste}, 
\ref{lem:alfabetaF-ste1} and 
\ref{lem:alfabetaF-ste2} that $\Gb$ is a Steinberg group, 
$\Gb\neq \,^{(3)}\!D_4(q)$; by the preceding Step, we also 
assume that $\Gb\neq \PSU_3(q)$.

\begin{stepa}\label{lem:alfabetaF-ste} 
If $q > 5$ and  
$\{\alpha,\beta\}\cap\{\theta(\alpha),\theta(\beta)\}=\emptyset$, 
then there exists $t_a\in \T\cap[\G^F,\G^F]$ 
for $a\in \I_4$ satisfying \eqref{eq:alfabetaF2}.
The same holds if $q=4$, except when $(\alpha,\theta(\alpha)) =-1$ and $(\theta(\alpha),\beta) = 1$.
\end{stepa}

Since $\theta$ preserves positivity of roots, 
we have $\theta(\alpha)\neq-\alpha$, $\theta(\beta)\neq-\beta$. 
Hence, for $m:=(\alpha,\theta(\alpha)), m':=(\theta(\alpha),\beta)$, we have $m,m'\in\{-1,0,1\}$.
Let $t_a$ be as in \eqref{eq:ta} with $\beta = \theta(\alpha)$. 
Then  $\alpha(t_a) = \zeta^{(a-1)(2+mq)}$, $\beta(t_a)=\zeta^{(a-1)(m'q-1)}$. 
Therefore,  \eqref{eq:alfabetaF2}, 
follows if $|\zeta^{(3+q(m-m'))}|\geq 4$. 
A direct estimate making use of the equalities 
$${\rm gcd}(q^2-1,q\pm3)={\rm gcd}(8,q\pm3), \mbox{ for }q\neq3$$
shows that  \eqref{eq:alfabetaF2} holds for $q > 5$, or for 
$q=4$ provided that $(m, m') \neq (-1, 1)$.

\begin{stepa}\label{lem:alfabetaF-ste1} 
If $q>7$ and $\beta\neq\theta(\alpha)$, then there exist 
$t_a\in \T\cap[\G^F,\G^F]$ for $a\in \I_4$ satisfying 
\eqref{eq:alfabetaF2}.
\end{stepa}

For $\alpha=\theta(\alpha)$, or $\alpha\not=\theta(\alpha)$ 
but $\beta=\theta(\beta)$, and $q>7$, 
the proof is as in Step \ref{lem:alfabetaF-chev}.
If $\alpha\not=\theta(\alpha)$ and $\beta\not=\theta(\beta)$, 
then this is Step \ref{lem:alfabetaF-ste}.

\begin{stepa}\label{lem:alfabetaF-ste2} 
If $q\not\in\{ 2,5,8\}$ and $\beta=\theta(\alpha)$, 
then there exist $t_a\in \T\cap[\G^F,\G^F]$ for $a\in \I_4$ 
satisfying \eqref{eq:alfabetaF2}.
\end{stepa}

This Step is proved as Step \ref{lem:alfabetaF-psu}. 

Hence the Proposition for a Steinberg group $\Gb$ 
different from $^{(3)}\!D_4(q)$ and $\PSU_3(q)$ 
follows from Steps \ref{lem:alfabetaF},  
\ref{lem:alfabetaF-ste}, \ref{lem:alfabetaF-ste1} and 
\ref{lem:alfabetaF-ste2}.

\begin{stepa}
Assume $G= ^{(3)}\!D_4(q)$ and $q\neq 2,3,4,7$. Then $\mathcal{O}$ is of
type F.
\end{stepa}

By the proof of Step 3 in Proposition \ref{prop:alfabetachevstein}
we can always assume that $\alpha=\alpha_2$ and $\beta\in\{\alpha_1,\alpha_1+\alpha_2+\alpha_3+\alpha_4\}$.
%Let us fix an ordering of the roots in increasing height ${\rm ht}$ and let 
%$u\in\mathcal{O}\cap \U^F$ be as in Definition \ref{def:alfabeta2}. 
%We consider the support of $u$ with respect to this ordering. 
%The outer automorphism $\theta$ of order $3$ permutes $\alpha_1,\alpha_3$
%and $\alpha_4$ and fixes $\alpha_2$. By inspection, for simple roots we have $\alpha\in \supp(u)$ 
%if and only if $\theta(\alpha)\in\supp(u)$. 
%In addition, $\gamma+\gamma'\not\in\Phi$ if ${\rm ht}(\gamma)={\rm ht}(\gamma')\geq2$ 
%or if ${\rm ht}(\gamma)={\rm ht}(\gamma')=1$ and $\gamma,\gamma'\neq\alpha_2$. 
%So, if $\{\alpha_1,\alpha_2\}\not\subset\supp (u)$ we have $\alpha \in \supp(u)$ 
%f and only if $\theta(\alpha) \in \supp (u)$ for every $\alpha\in\Phi^+$. 
%Thus, if $\alpha_2\not\in\supp(u)$, condition (\ref{eq:ab-property}) is not verified 
%for any pair $\alpha,\beta\in\supp(u)$ such that $\alpha+\beta\in \Phi$. 
%So, $\alpha_2\in\supp(u)$. 
%If $\alpha_1\in\supp(u)$ then we take $\alpha=\alpha_2$, $\beta=\alpha_1$. 
Let $\zeta \in \mathbb{F}_{q}^{\times}$ be a generator. 
Let $e_a=a-1$, with $a\in\I_4$. In the first case we 
take $t_a=\beta^{\vee}(\zeta)\theta(\beta)^{\vee}(\zeta)\theta^{2}(\beta)^{\vee}(\zeta)
\alpha^{\vee}(\zeta^{e_a}) $.  
Since $\zeta^q=\zeta$, we have $t_a \in \mathbb{T}^F$. Further,
since $(\alpha,\theta^{i}(\beta)) = -1 $ and 
$(\theta^{i}(\beta),\theta^{i}(\beta))=2$
we have
$\alpha(t_a)\beta(t_b) = \zeta^{-1+2e_a-e_b}$.
%  and
%$\alpha(t_b)\beta(t_a)=\zeta^{-3+2e_b}\zeta^{2-e_a}$ and
%$(\alpha(t_a),\beta(t_a)) = ( \zeta^{-3+2e_a},\zeta^{2-e_a}) $. 
As in Step 2, \eqref{eq:alfabetaF2} 
are satisfied if $q\neq 2,3,4,7$.
% If, instead, $\alpha_1\not\in\supp(u)$, then $u$ 
%has the $\alpha\beta$-property if and only if $\alpha_1+\alpha_2+\alpha_3+\alpha_4\in\supp(u)$ 
%and $\{\alpha_1+\alpha_2,\alpha_2+\alpha_3+\alpha_4\}\not\subset\supp(u)$. 
In the second case we take $\alpha=\alpha_2$, $\beta=\alpha_1+\alpha_2+\alpha_3+\alpha_4$.  
Then, the proof follows as in Step 2. 
Indeed, define $t_a=\alpha^\vee(\zeta^{e_a}) \in\T^F\cap [\G^F,\G^F]$.
Then $\alpha(t_a)\beta(t_b)=\zeta^{2e_a-e_b}$  and a
direct verification gives the claim.
Thus, $\mathcal{O}$ is of
type F by Step 1.
\epf

\subsection{Regular unipotent classes}\label{subsec:preliminaries-regular}

A \emph{regular} unipotent conjugacy class in a reductive 
algebraic group is the unique unipotent class with maximal dimension.
Then we say that $\Oc$ is \emph{regular} if it is contained  in the regular unipotent class in $\G$.
We shall prove often that a class is of type D or F
by considering the intersection with a smaller group, 
containing a regular class of the latter. Thus, we need
to see when regular classes are of type D or F.

\begin{prop}\label{prop:regular} Assume that $q$ is odd. 
Let $\Gb$ be a finite simple group of Lie type 
not of type $A_1$. If $\Oc$ is regular, then it is of type D.
\end{prop}

\pf By \cite[3.2,\,3.3]{ihes},
every regular unipotent element $u$ in $\U$ can be written as
$u = \U' x_{\alpha_1}(a_1)\cdots x_{\alpha_n}(a_n)$ 
where $\U'$ is the product of root subgroups
of height at least $2$ and each $a_i\in\kk^\times$.
If the rank  of $\G$ is not $1$, this ensures that
for every $u\in \Oc\cap\U^{F}$, there are $\alpha$, 
$\beta$ simple adjacent roots in  $\supp(u)$;
hence $\Oc$ has
the $\alpha\beta$-property. Now Proposition \ref{prop:alfabetachevstein} applies.
\epf

\begin{prop}\label{prop:regular-char2} Assume $\Gb$ is either
\begin{enumerate}
  \renewcommand{\theenumi}{\alph{enumi}}\renewcommand{\labelenumi}{(\theenumi)}
\item\label{item:regular-char2-chevalley} a Chevalley group with $q>7$ and $\G\neq\SL_2(\kk)$, or 
\item\label{item:regular-char2-psu3} $\PSU_3(q)$, with $q\not\in\{2,5,8\}$, or
\item\label{item:regular-char2-psun} $\PSU_n(q)$, with $n\geq5$, or 
$\,^{(2)}\!E_6(q)$, and $q\not\in\{2,3,5\}$, or
\item\label{item:regular-char2-psu4-dn} $\,^{(2)}\!D_n(q)$ for 
$n\geq4$ or $\PSU_4(q)$, and $q>7$,  or 
\item $\,^{(3)}\!D_4(q)$ and $q\neq 2,3,4,7$.
\end{enumerate}
Then every regular unipotent class in $\Gb$ is of type F.
\end{prop}
\pf  Arguing as in Proposition \ref{prop:regular}
there exists $u\in \Oc\cap\U^{F}$ and $\alpha,\beta \in \Phi^{+}$
such that $\Oc$ has
the $\alpha\beta$-property, hence we may invoke 
Proposition \ref{prop:alfabetaF}. If $\Gb=\PSU_n(q)$, 
with $n\geq5$, or $\Gb=\,^{(2)}\!E_6(q)$
we can always find adjacent simple roots $\alpha$ 
and $\beta$ such that $\{\alpha,\beta\}\cap \{\theta(\alpha),\theta(\beta)\}=\emptyset$, 
so Step \ref{lem:alfabetaF-ste}  applies. If $\Gb=\,^{(2)}\!D_n(q)$ 
for $n\geq4$ or $\PSU_4(q)$, 
then we can always find adjacent simple roots $\alpha$
and $\beta$ with $\beta \neq \theta(\alpha)$ and  Step \ref{lem:alfabetaF-ste1}  applies.
\epf

\begin{obs}If $p$ is good (see \cite[I.4.3]{sp-st} 
for the list of bad primes) then all regular unipotent classes in $\G^F$ 
are isomorphic as racks \cite[Lemma 4.1]{TZ}. 
But this is not always the case for $p$ bad. Let, for instance, $p=2$, 
 $\G^F$ $=\Sp_4(2)\cong \Sim_6$. 
The regular unipotent classes $\Oc$ and $\Oc'$ 
correspond to the partitions $(1^2, 4)$ and $(2, 4)$, 
and have isomorphic centralizers. We compute the 
inner groups of them, see \cite[Definition 1.3]{AG-adv}, using
\cite[Lemma 1.9]{AG-adv}:  $\Inn_{\Oc} =\Sim_6$ and 
$\Inn_{\Oc'} =\as$. Thus $\Oc$ and $\Oc'$ are not isomorphic as racks.
\end{obs}

We next deal with some specific groups. 
See \S \ref{sec:symplectic} for the needed notation of symplectic groups. 

\begin{lema}\label{lem:sp-regular} Let $q>2$ be even. 
The regular unipotent classes in $\Sp_{2n}(q)$ are of type F.
\end{lema}
\pf  There are exactly $2$ regular unipotent classes in 
$\Sp_{2n}(q)$  \cite[Theorem 6.2.1]{LS}. 
Both are treated similarly, so fix one of them, say $\C$. 
There is an upper-triangular matrix $u\in \C$. 
By Jordan theory $u$ is  regular in $\SL_{2n}(q)$ \cite[IV.2.15.9(ii)]{sp-st}, 
so all its coefficients in the upper subdiagonal are $\neq 0$. 

Assume first that $n=2$. Then we may assume that 
$u = \left(\begin{smallmatrix}
      1 & x & 0 & p\\
      0 & 1 & y & xy\\
      0 & 0 & 1 & x\\
      0 & 0 & 0 & 1
     \end{smallmatrix}\right)
$. Indeed, if $u = \left(\begin{smallmatrix}
      1 & x & l & p\\
      0 & 1 & y & m\\
      0 & 0 & 1 & z\\
      0 & 0 & 0 & 1
     \end{smallmatrix}\right)
$, then since $u\in \Sp_4(q)$ and it is a regular element we have that
$x=z$, $m=l+xy$ and $xyz\neq 0$. Conjugating by 
$v = \left(\begin{smallmatrix}
      1 & ly^{-1} & 0 & 0\\
      0 & 1 & 0 & 0\\
      0 & 0 & 1 & ly^{-1}\\
      0 & 0 & 0 & 1
     \end{smallmatrix}\right)$ we obtain that $(vuv^{-1})_{13}= 0$
     and $(vuv^{-1})_{24} = xy \neq 0$.
Let $\zeta$ be a generator of $\F_q^\times$.
Conjugation by the diagonal matrices 
$(\zeta^a,\zeta^b,\zeta^{-b},\zeta^{-a})$, $a,b\in\{0,1\}$,
provides $4$ different representatives $x_{a,b}$ of $\Oc_x^{\Sp_{4}(q)}$.  
We check that 
the subracks $R_{a,b}:=\U^F\trid x_{a,b}$ are disjoint 
for $(a,b)\neq(c,d)$ as in \cite[3.1]{ACG-I}. A direct computation looking at 
$(x_{a,b}x_{c,d})_{13}$ and  $(x_{a,b}x_{c,d})_{14}$ 
shows that $x_{a,b}\trid x_{c,d}\neq x_{c,d}$ if $(a,b)\neq(c,d)$ for all $q>2$.

Assume now $n>2$. Then Lemma \ref{lem:parabolic} applies with $P = \Pa^{F}$ such that
$\Pa$ is the standard
parabolic subgroup containing $\mathbb{L} = \T\Sp_{4}(\Bbbk)$ as a Levi factor and 
$L = \mathbb{L}^{F}$, the $\F_{q}$-points of $\mathbb{L}$. Indeed since $u$ is regular unipotent, 
it is a $p$-element so its image 
$\overline{u}$ is contained in $\Sp_4(q)\subset L$ and it is regular therein.
\epf

For further use, we treat here regular classes in other groups. 
Recall the notation in \S \ref{subsubsec:unitary}.

\begin{lema}\label{lem:GU} Let $q=2^{2h+1}$, where $h\in \N_0$. 
The regular unipotent classes in $\GU_n(q)$ for $1<n$ odd are of type D.
\end{lema}

\pf Assume $q=2$, $n=3$. By Subsection \ref{subsubsec:unitary}, 
there is a unique unipotent 
regular conjugacy class $\Oc$ in $\GU_3(2)$. 
Let $\zeta,\,\eta\in \F_4^{\times} - 1$. Then
$r=\left(\begin{smallmatrix}
1&1&\zeta\\
0&1&1\\
0&0&1
\end{smallmatrix}\right)\in\Oc$. 
By \cite[6.22]{Hu}, $C_{\SL_3(\kk)}(u)=Z(\SL_3(\kk))C_{\SL_3(\kk)}(u)^\circ$. 
It is not hard to verify that, for $F$ the twisted 
Steinberg endomorphism on $\SL_3(\kk)$, the $F$-twisted 
action of $Z(\SL_3(\kk))\simeq \Z/3$ on itself is trivial, see \cite[Example 21.14]{MT}. 
Thus, there are exactly $3$ regular unipotent conjugacy 
classes in $\SU_3(2)$. Let $x\in \kk$ be such that $x^3=\eta^{-1}$. 
Then  for $g=\left(\begin{smallmatrix}
x^4\\
&x\\
&&x^4\end{smallmatrix}\right)$ we have $g^{-1}F(g)=\eta\id_3$ 
so, $grg^{-1}\in \Oc\setminus \Oc_r^{\SU_3(2)}$.
Clearly ${\Jf}_3\in \SU_3(2)<\GU_3(2)$, so also  
$s={\Jf}_3 g r g^{-1} {\Jf}_3\in \Oc\setminus \Oc_r^{\SU_3(2)}$.
By a direct computation, $(rs)^2\neq (sr)^2$. 
Since $\langle r,s\rangle\subset \SU_3(2)$,  $\Oc$ is of type D.

% while \eqref{item:psl3-regular}
%follows since $\SU_3(2) < \SL_3(4) < \SL_3(2^{2m})$.

Assume $q=2$, $n=2l+1>3$. Let $\Pa$ be the standard $F$-stable 
parabolic subgroup associated to the simple roots $\alpha_l,\alpha_{l+1}$ 
and let $\Le$ be the corresponding standard $F$-stable Levi subgroup; 
$\Le$ contains a subgroup isomorphic to $\GL_3(\kk)$. 
Then, $\Le^F$ contains a subgroup isomorphic to $\GU_3(q)$ and  Lemma \ref{lem:parabolic} applies with $P=\Pa^F$ and $L=\Le^F$, by the  case $n=3$.
The claim for general $q$ follows since $\GU_n(2)< \GU_n(2^{2h+1})$.
\epf

\subsection{Further remarks}\label{subsec:preliminaries-remarks}
We shall often invoke the following result. 

\begin{lema}\label{lem:SL}\cite[\S 3.5]{ACG-I} Let $\Oc$ be a 
unipotent class in $\SL_n(q)$,  with partition $(\lambda_1,\ldots,\lambda_k)$. 
\ref{tab:unipotent-chevalley-collapse-typeA} summarizes 
when $\Oc$ is of type D or F. \qed
\end{lema}

\begin{align}\label{tab:unipotent-chevalley-collapse-typeA}
\stepcounter{tabla}\tag*{Table \thetabla}
\begin{tabular}{|c|c|c|c|}
\hline $n$   & $q$ & type $(\lambda_1, \dots, \lambda_k)$ & Type  \\
\hline
$2$   & odd square $> 9$      & $(2)$ &  D \\
 \hline\hline
$ > 2$  & odd   & $\lambda_{1} \ge 3$ &  D   \\
  \cline{3-4}
   &    &  $(2, 2, \dots)$ &  D   \\
  \cline{3-4}
   &    &  $(2, 1 \dots)$ &  D   \\
 \hline\hline
$> 2$ &  even & $\lambda_1 \ge 5$ & F  \\
 \cline{3-4}
 &   & $\lambda_1 = 4$ & D 
\\ \cline{3-4}
  &     & $(3, 3, \dots)$  & D \\ 
  \cline{3-4}
  &     & $(3, 2, \dots)$  & F \\
 \cline{3-4}
 &   & $(3, 1, \dots)$ & D
\\ \cline{3-4}
 &   & $(2, 2, \dots)$ & D\\ 
\cline{3-4}
 &   & $(2,1,1,1,\dots)$ &  F
\\
\hline\hline
$3$ & even $\geq 8$ & $(3)$ &   F\\
\cline{2-4}
  & 4 & $(3)$  &D
\\
%\cline{4-4}\cline{2-2}
 \hline
 \end{tabular}
 \end{align}

We end this subsection with another useful observation.

\begin{obs}\label{rem:parabolic} Let $\Pa$ be an $F$-stable 
parabolic subgroup of $\G$, let $\Le$ be an $F$-stable Levi subgroup  
and let $\pi\colon\Pa\to\Le$ be the  projection associated with 
the Levi decomposition $\Pa = \Le\urad$. 
Let $G = \G^F$, $P = \Pa^F$, $Q = \urad^F$, $L = \Le^F$.
\begin{enumerate}
\item Let $r,s\in P$ with $s\in Q\not\ni r$. Then 
$\Oc_r^{\langle r,s\rangle}\neq \Oc_s^{\langle r,s\rangle}$ because $Q\nsbgp P$. 
\item If moreover $\Oc_r^{G}=\Oc_s^{G}$ and 
$\pi(rs)^2\neq \pi(sr)^2$, then $\Oc_r^{G}$ is of type D. 
\end{enumerate}
\end{obs}

\section{Unipotent classes in  finite symplectic groups}\label{sec:symplectic}

In this section, $\G$ is the symplectic group $\Sp_{2n}(\kk)$, 
that is the subgroup of $\GL_{2n}(\kk)$ leaving invariant the
bilinear form  $\left(\begin{smallmatrix}
0&{\Jf}_n\\
-{\Jf}_n&0\end{smallmatrix}\right)$, for  ${\Jf}_n=
\left(\begin{smallmatrix}
&&1\\
&\iddots&\\
1&&
\end{smallmatrix}\right)$. We assume $n\geq2$, since $\Sp_2(\kk)=\SL_2(\kk)$.
Let $\B$ be the Borel subgroup of $\G$ consisting of 
upper triangular matrices.
Since $\G$ is simply connected,  $\G^F=[\G^F,\G^F] = \Sp_{2n}(q) =: G$ 
and $\Gb= \PSp_{2n}(q) = \G^F/Z(\G^F)$. 
By the isogeny argument, Lemma \ref{lem:isogeny-argument},
it suffices to consider unipotent classes in $G$.

\subsection{Symplectic groups for $q$ odd}\label{subsec:odd-symplectic}

To a unipotent class $\Oc$ in $\Sp_{2n}(\kk)$ we attach the partition of $2n$ 
determined by the Jordan form of $\Oc$
in $\GL_{2n}(\kk)$. If $q$ is odd, then the partition uniquely determines $\Oc$. 
The partitions corresponding to unipotent classes in $\G$ are of the form
$(1^{r_1},2^{r_2},\ldots,2n^{r_{2n}})$ where $r_i$ is even for every 
odd $i$  \cite[IV.2.15.9(ii)]{sp-st}. We call them \emph{symplectic} partitions.

Let $u\in \G$ unipotent. There is a reductive subgroup $\J$ of $\G$ containing $u$ 
as a regular unipotent element, such that $C_{\G}(u)=C_{\G}(\J)\V$ 
where $\V$ is a connected normal subgroup of $C_{\G}(u)$ \cite[Lemmata 3.14, 3.17]{LS}. 
Namely, $\J$ is given in \cite[(3.4), p. 48]{LS}: if
$\Oc_u$ corresponds to $(1^{r_1},2^{r_2},\ldots,2n^{r_{2n}})$, then
\begin{align}\label{eq:def-J}
\J\cong \prod_{i\mbox{ odd }}{\mathbf{O}}_i(\kk)\times \prod_{i\mbox{ even }}\Sp_i(\kk),
\end{align}
where the product is taken over those $i$ such that $r_{i}\neq 0$.
We can always assume that $\J$ is $F$-stable and that $F$ induces an $\F_q$-split 
morphism on each of its simple factors \cite[p. 113]{LS}.
Recall that $\clas{G}{u}$ denotes the set of $G$-conjugacy classes 
contained in $\Oc_u^{\G}$, when $u\in G$.
%; see the discussion in \S \ref{subsubsec:preliminaries-classes}.

\begin{lema}\label{lem:sympl-odd} Let $u$ be a nontrivial unipotent element in $G$ 
associated with the partition $(1^{r_1},\ldots,\,n^{r_n})$.  
Assume that one of these conditions hold:
\begin{enumerate}
\item\label{item:1} there exists $i>3$ for which $r_i\neq0$;
\item\label{item:4} $9 < q$ is a square and the partition is $(1^{r_1},2^{r_2})$ with $r_2>0$.
\end{enumerate}
 Then $\Oc$ is of type D for every  $\Oc\in \clas{G}{u}$.
 \end{lema}
\pf Since $u$ is unipotent, it lies in the following subgroup of $\J$
$$\M=\prod_{i\mbox{ odd }}{\mathbf{SO}}_i(\kk)\times \prod_{i\mbox{ even }}\Sp_i(\kk)$$ and
each component of $u$ in $\M$ is regular in its  factor. We show that $\Oc_u^{\M^F}$ is of type D. 
Case \eqref{item:1} follows from Proposition~\ref{prop:regular}.
In Case \eqref{item:4}, $r_2>0$ and  $\M^F$ is a group of type $A_1$; 
hence \cite[Lemma 3.6]{ACG-I} applies.
For the other classes in $\clas{G}{u}$, we apply Lemma \ref{lem:J} $(c)$; indeed
$C_{\G}(u)=C_{\G}(\J)\V$ and $\V$ is connected, so representatives of 
$A(u)$ can be found in $C_{\G}(\J) < C_{\G}(\M)$.
\epf

By Lemma \ref{lem:sympl-odd} \eqref{item:1},  it remains to consider 
the partitions $(1^{r_1},2^{r_2}, 3^{r_3})$.
We start by $(1^{r_1}, 3^{r_3})$; the argument in 
Lemma \ref{lem:sympl-odd} \eqref{item:4} also applies for it, but
there is an alternative without the restrictions in the parameters.

\begin{lema}\label{lem:sp31} Let $u\in G$ be a unipotent element 
corresponding to a partition of the form $(1^{r_1}, 3^{r_3})$, with $r_3 > 0$. 
Then $\Oc_u^{G}$ is of type D.
\end{lema}
\pf By \cite[Theorem 3.1 (v)]{LS} we have
$C_\G(u)=C_{\G}(u)^\circ$ so  $\clas{G}{u}$ consists of a single class $\Oc$. We set $r_3=2e$.
Let $v=\left(\begin{smallmatrix}
1&1&0\\
0&1&1\\
0&0&1
\end{smallmatrix}\right)$. A representative of $\Oc$ is
$\left(\begin{smallmatrix}
v\otimes \id_e&0&0\\
0&\id_{2n - 3r_3}&0\\
0&0&({\Jf}_3\,^t\!v^{-1}{\Jf}_3)\otimes \id_e\end{smallmatrix}\right)$. 
Now we may assume that $e=1$; in this case
the injective morphism $\iota: \SL_3(q)\to G$,
$X\mapsto \diag(X, \id_{2n-6},{\Jf}_3\,^t\!X^{-1}{\Jf}_3)$
induces a rack embedding $\Oc_v^{\SL_3(q)}\hookrightarrow \Oc_u^{G}$ 
and by the isogeny argument, Proposition~\ref{prop:regular} applies.
\epf

\begin{lema}\label{lem:sp-2columns}
Let $\Oc$ be a unipotent class in $G$ with partition $(1^{r_1}, 2^{r_2})$, with $r_2>1$.
If $q>3$ or  $n>2$, then the class is of type D.
\end{lema}
\pf In this case $A(u)\simeq \Z/2$ \cite[Theorem 3.1 (v)]{LS}, 
so $\clas{G}{u}$ consists of 2 classes.
One of them is represented by
$u=\left(\begin{smallmatrix}\id_{r_2}&0&{\Jf}_{r_2}\cr
0&\id_{r_1}&0\\
0&0&\id_{r_2}\end{smallmatrix}\right)$. We find a  representative for the other.
Recall the notation \eqref{eq:def-J}; it can be shown that 
$C_{\G}(\J)\simeq {\mathbf{O}}_{r_2}(\kk)\times \Sp_{r_1}(\kk)$ 
is the subgroup of matrices
$g_{A,M}=\left(\begin{smallmatrix}A&0&0\\
0&M&0\\
0&0&{\Jf}_{r_2}A{\Jf}_{r_2}\end{smallmatrix}\right)$ with 
$A^t\!A=\id_{r_2}$ and $M\in\Sp_{r_1}(\kk)$.
The nontrivial element in $A(u)$ is represented by
$g_{L,\id_{r_1}}$ for $L=\diag(-1, 1, \dots, 1)$. 
Let $\xi\in\F_{q^2}\setminus\F_{q}$ be such that $\xi^{q-1}=-1$, 
so $\zeta=\xi^2\in\F_q$ is not a square in $\F_q$. 
Let $g=\diag(\xi,1, ,\dots, 1,\xi^{-1})\in\Sp_{2n}(\kk)$. 
Then $g^{-1}F(g)=g_{L,\id_{r_1}}$, so
$\clas{G}{u}=\{\Oc_u^{G},\Oc_{v}^{G}\}$ where
$$v=gug^{-1}=
\left(\begin{smallmatrix}1&&&&\zeta\\
&\id_{r_2-1}&&{\Jf}_{r_2-1}\\
&&\id_{r_1}\\
&&&\id_{r_2-1}\\
&&&&1\end{smallmatrix}\right).
$$

Let $\M$ be the $F$-stable subgroup of $\Sp_{2n}(\kk)$ of matrices 
$\left(\begin{smallmatrix}a&0&b\\
0&M&0\\
c&0&d\\
\end{smallmatrix}\right)$ with $ad-bc=1$ and $M\in\Sp_{2n-2}(\kk)$. 
Clearly, $\M\simeq \SL_2(\kk)\times \Sp_{2n-2}(\kk)$. 
Then $\Oc_u^{\M^F}\simeq\Oc_{v}^{\M^F}\simeq\Oc_{x_1}^{\SL_2(q)}\times \Oc_{y_1}^{\Sp_{2n-2}(q)}$ with
$$x_1=\left(\begin{smallmatrix} 1&1\\ 0&1\end{smallmatrix}\right),\quad y_1
=\left(\begin{smallmatrix}\id_{r_2-1}&0&{\Jf}_{r_2-1}\cr
0&\id_{r_1}&0\\
0&0&\id_{r_2-1}\end{smallmatrix}\right).$$ We show that this subrack is of 
type D by application of Lemma \ref{lema:ACG-2-10}.
First, $x_2=\left(\begin{smallmatrix} 1&0\\ -1&1\end{smallmatrix}\right)\in \Oc_{x_1}^{\SL_2(q)}$ 
satisfies $(x_1x_2)^2\neq (x_2x_1)^2$.
So, we have to find  $y_2\in \Oc_{y_1}^{\Sp_{2n-2}(q)}$ commuting with $y_1$.

Assume $q > 3$, and let $\eta\in \F_q^\times\setminus \{1,-1\}$. Then we take
$$y_2=\left(\begin{smallmatrix}\eta\id_{r_2-1}\\
&\id_{r_1}\\
&&\eta^{-1}\id_{r_2-1}\end{smallmatrix}\right)\trid y_1= 
\left(\begin{smallmatrix}\id_{r_2-1}&0&{\eta^2\Jf}_{r_2-1}\cr
0&\id_{r_1}&0\\
0&0&\id_{r_2-1}\end{smallmatrix}\right).$$ 
Assume now  that $n>2$. If $r_2>2$, then
$$y_2=\left(\begin{smallmatrix}
1&1\\
&1\\
&&\id_{2n-6}\\
&&&1&-1\\
&&&&1
\end{smallmatrix}\right)\trid y_1\in \Oc_{y_1}^{\Sp_{2n-2}(q)}$$ commutes with $y_1$.
If $r_2=2$, then necessarily $r_1>1$.
In this case we take
$y_2=\left(\begin{smallmatrix}
0&1\\
-1&0\\
&&\id_{2n-6}\\
&&&0&-1\\
&&&1&0
\end{smallmatrix}\right)\trid y_1=\left(\begin{smallmatrix}1&&&0&0\\
&1&&1&0\\
&&\id_{2n-6}\\
&&&1\\
&&&&1\end{smallmatrix}\right).$
\epf

\begin{lema}\label{lem:321sp} 
Let $u$ be a unipotent element in $G$ with partition 
$(1^{r_1}, 2^{r_2}, 3^{r_3})$, such that $r_2r_3 > 0$. 
Then $\Oc_u^{G}$ is of type D.
\end{lema}
\pf Here $r_3=2a$ is even and $C_{\G}(\J)\simeq \Sp_{r_1}(\kk)
\times \Ort_{r_2}(\kk)\times\Sp_{r_3}(\kk)$, so 
$A(u)\simeq\Z/2$ \cite[Theorem 3.1(v)]{LS}
and $\clas{G}{u} =\{\Oc_u^{G},\Oc_{v}^{G}\}$ has 2 elements. Let
$x=\left(\begin{smallmatrix}1&1&0\\
0&1&1\\
0&0&1\\
\end{smallmatrix}\right)$ and $w$ in $\Sp_{2n-3r_3}(q)$ 
unipotent with partition $(1^{r_1}, 2^{r_2})$. Then we choose
$$u=\left(\begin{smallmatrix}
x\otimes \id_a\\
&w\\
&&{\Jf_{3}}^t\!x^{-1}{\Jf_{3}}\otimes \id_{a}\\
\end{smallmatrix}\right).
$$
Let $\M < \G$ consisting
of matrices of the form $\diag(\id_{3a},Y,\id_{3a})$, with $Y$ in $\Sp_{2n-3r_3}(\kk)$, and 
let $\J_w$ be  the reductive group containing $w$ as a regular element, cf. \cite[3.2.3]{LS}.
Then $C_{\M}(\J_w)\simeq \Sp_{r_1}(\kk)\times \Ort_{r_2}(\kk) < C_{\G}(\J) < C_{\G}(u)$ and 
a representative $b$ for the nontrivial class in $A(u)$ may
be chosen in $\M$. By Lang-Steinberg's Theorem, there exists $g\in \M$ 
such that $g^{-1}F(g)=b$, so we may pick $v = gug^{-1}$. In other terms, 
$v=\left(\begin{smallmatrix}
x\otimes \id_a\\
&z\\
&&{\Jf}_{3}^t\!x^{-1}{\Jf}_{3}\otimes \id_a\\
\end{smallmatrix}\right)$ with $\clas{\Sp_{2n-3r_3}(\kk)}{w} 
= \big\{\Oc_{w}^{\Sp_{2n-3r_3}(q)},\Oc_{z}^{\Sp_{2n-3r_3}(q)}\big\}$.
Then $\Oc_u^{G}$ and $\Oc_{v}^{G}$ contain a subrack isomorphic to 
$\Oc_{x}^{\SL_3(q)}$, and 
by the isogeny argument, Proposition \ref{prop:regular} applies.
\epf

\begin{lema}\label{lem:sp43}There are two unipotent classes in 
$\G^{F}=\Sp_4(3)$ of type $(2^{2})$;
one of them is of type D and the other is cthulhu.
\end{lema}
\pf As in the proof of Lemma \ref{lem:sp-2columns}, 
there are two classes of type $(2^2)$, represented by 
$z=\left(\begin{smallmatrix}1&&&1\\
&1&1\cr
&&1\\
&&&1\\
\end{smallmatrix}\right)$ and $w=\left(\begin{smallmatrix}1&&&-1\\
&1&1\cr
&&1\\
&&&1\\
\end{smallmatrix}\right)$.
It can be verified that 
$x=x_{\alpha_1}(1)=\left(\begin{smallmatrix}1&1\\
&1&\cr
&&1&-1\\
&&&1\\
\end{smallmatrix}\right)\in\Oc_w^{\Sp_4(3)}$. The discussion in \cite[3.1]{ACG-I}
 shows that $\Oc_x^{\U^F}\not=
\Oc_w^{\U^F}$, and that $xw\not=wx$. In addition, 
$xw$ and $w(xw)w^{-1}=wx$ lie in $\U^F$ so they have odd order, hence $(xw)^2\neq (wx)^2$. 
The claim on the other class was verified with GAP.
\epf

\subsection{Symplectic groups for $q$ even}\label{subsec:even-symplectic}

In this section $q$ is even, so the symplectic group $\G$ is the subgroup 
of $\GL_{2n}(\kk)$ leaving invariant the
bilinear form  $\Jf_{2n}$. 
Here symplectic partitions do not distinguish conjugacy classes.

\subsubsection{Unipotent conjugacy classes}\label{subsubsec:even-symplectic-classes}
We parametrize the unipotent conjugacy classes in $\G$ as in \cite[6.1, cf. Lemma 6.2]{LS}. 
Let $V$ be the natural representation of $\G$ and let $u\in \G$ unipotent. 
Then $V$ decomposes, as an $u$-module by restriction, 
into an orthogonal direct sum of indecomposable submodules 
(where $k, r \in \N_0$, the $m_i$'s are distinct, ditto for the $k_j$'s)
\begin{align}\label{deco}
V &= \bigoplus_{i=1}^k W(m_i)^{a_i}\oplus \bigoplus_{j=1}^r V(2k_j)^{b_j},&  &0 <a_i,\; 0<b_j\leq 2.
\end{align}
We describe the summands in the right-hand side:

\noindent $\bullet \ \dim W(m_i) = 2m_i$ and $u_{\vert W(m_i)}$ is regular  in a subgroup $\Hb_{m_i}$, that
is the image of $\SL_{m_i}(\kk)$  by the embedding of $\GL_{m_i}(\kk)$ in $\Sp(W(m_i))$ given by 
\begin{align}\label{eq:image-embedding}
X \mapsto  \diag(X, {\Jf}_{m_i}\,^t\!X^{-1}{\Jf}_{m_i})
\end{align}
and thus $u_{\vert W(m_i)}$ is of partition $(m_{i},m_{i})$ in $\Sp(W(m_i))$;

\smallbreak
\noindent $\bullet \ \dim V(2k_j) = 2k_j$ and $u_{\vert V(2k_j)}$ 
is regular  in a  subgroup $\J_{2k_j}\simeq\Sp_{2k_j}(\kk)$
and thus $u_{\vert V(2k_j)}$ is  of partition $(2k_j)$. 

\smallbreak
Set  $\Hb = \prod_{i=1}^k \Hb_{m_i}^{a_i}$, $\J=\prod_{j=1}^r \J_{2k_j}^{b_j}$. 
Then $u$ is regular in $\M:= \Hb \times \J$.
Let $\W=\bigoplus_{i=1}^k W(m_i)^{a_i}$ and $\Vc=\bigoplus_{j=1}^r V(2k_j)^{b_j}$. Then 
\begin{align}\label{eq:cadena-M}
 \M <
\prod_{i=1}^k\Sp(W(m_i)^{a_i})\times \prod_{j=1}^r\Sp(V(2k_j)^{b_j}) <
\Sp(\W)\times \Sp(\Vc) < \G.
\end{align}
 
By the description in \cite[p. 91]{LS}, there is $u\in G = \G^F$ such that all subgroups 
$\Hb_{m_i}$,  $\J_{2k_j}$, $\Hb$, $\J$, $\M$,  $\Sp(\W)$,  $\Sp(\Vc)$
are $F$-stable and $F$  acts on each of them by a split Frobenius automorphism. In particular,
\begin{align*}
\Hb^F &\simeq \prod_{i=1}^k \SL_{m_i}(q)^{a_i}, &\J^F &\simeq \prod_{j=1}^r \Sp_{2k_j}(q)^{b_j}.
\end{align*}
We fix this $u$ in the rest of this Subsection.

\subsubsection{Representatives of classes in $\clas{G}{u}$}
\label{subsubsec:even-symplectic-A(u)}

We now address the problem of finding suitable subracks for $\Oc\in\clas{G}{u}$,
that we recall again is the set of $G$-conjugacy classes 
contained in $\Oc_u^{\G}$.
First we need some information on $A(u)$, cf. \eqref{eq:a(x)}.

\begin{lema}\label{lema:representatives-a(u)-evensp}
There is a set of representatives $\Xi$ of $A(u)$ in $C_{\G}(u)$ such that for 
every $x\in \Xi$ there is  $g\in \G$ with $x = g^{-1}F(g)$ satisfying: 
\begin{enumerate}\renewcommand{\theenumi}{\alph{enumi}}
\renewcommand{\labelenumi}{(\theenumi)}
\item\label{item:f-stable} $F(g\M g^{-1})=g\M g^{-1}$ 
and $F(g\J_{2k_j} g^{-1})=g\J_{2k_j} g^{-1}$ for every $j$.

\smallbreak
\item\label{item:invaraints} $(g\M g^{-1})^F\simeq  
\prod_{i=1}^k (\SL_{m_i}(q)^{a_i-1}\times G_i) \times \J^F$,
where $G_i$ is  $\SL_{m_i}(q)$ or $\SU_{m_i}(q)$;  $\SU_{m_i}(q)$ 
may occur only if $m_i>1$  is odd.
\end{enumerate}
\end{lema}
\pf  By the proof of \cite[Theorem 6.21]{LS}, there is a 
maximal torus $\T_0$ of $C_{\G}(u)$ such that
\begin{align}\label{eq:even-symplectic-A(u)}
 C_{\G}(u)=C_{\G}(u)^\circ N H
\end{align}
where   $N=N_{\G}(\T_0)\cap C_{\Sp(\W)}(u)$, $H=C_{\Sp(\Vc)}(u)$, 
$NH \simeq H\times N$. Also, 
\begin{align*}
(N H)\cap C_{\G}(u)^\circ &= (N\cap C_{\G}(u)^\circ)  (H\cap C_{\G}(u)^\circ), 
& H/H\cap C_{\G}(u)^\circ &= H/H^\circ.
\end{align*}
Hence, we may construct a set of representatives $\Xi$ for $A(u)$  
as a product $\Sigma\Sigma'\in NH$ where $\Sigma$, resp. $\Sigma'$, is a set  of
representatives of $N/N\cap C_{\G}(u)^\circ$, respectively of  $H/H^\circ$.
We claim that  there are $\Sigma,\Sigma'\subset N_{\G}(\M)\cap N_{\G}(\J_{2k_j})$ for every $j$.
First $H/H^\circ$ is generated by images of the components $u_j$ of $u$ in some 
factors $\J_{2k_j}$ \cite[Lemmata 6.13, 6.14]{LS}, so any $\Sigma'\subset\J$ will do. 
For $\Sigma$ we need additional information from the proof of \cite[Theorem 6.21]{LS}:
\begin{itemize}
\item $\T_0$ is a product of subtori $\T_i$ of dimension $a_i$ acting on a 
single summand $W(m_i)^{a_i}$ without fixed points;
\item $N=\prod_{i=1}^k N_i$ where $N_i=N_{\G}(\T_i)\cap C_{\Sp(W(m_i)^{a_i})}(u)$.
\end{itemize}

Then $\prod_{i=1}^k N_i/(N_i\cap C_{\G}(u)^\circ)$ maps onto 
$N/(N\cap C_{\G}(u)^\circ)$ and  $N_i\subset N_{\G}(\Sp(W(m_j)^{a_j})
\cap N_{\G}(\J_{2k_j})$ for every $i, j$.
In order to describe the action of an $x\in N_i$ on $\Sp(W(m_i)^{a_i})$ 
we analyze the component of $u$ lying in this subgroup. 
We may assume it is $\left(\begin{smallmatrix} 
v\otimes \id_{a_i}\\
& {\Jf}_{m_{i}}\;^{t}v^{-1}{\Jf}_{m_{i}}\otimes \id_{a_i}\end{smallmatrix}\right)$
where $v$ is a regular unipotent matrix in $\SL_{m_i}(\kk)$,
and so  $\T_i$ is the subgroup of diagonal matrices 
$$(\lambda_1\id_{m_i},\,\ldots,\,\lambda_{a_i} \id_{m_i},
\lambda_{a_i}^{-1}\id_{m_i},\,\ldots,\,\lambda_{1}^{-1} \id_{m_i}).$$

Then $N_{\Sp_{2m_ia_i}(\kk)}(\T_i)$ normalizes 
$[C_{\Sp_{2m_ia_i}(\kk)}(\T_i),C_{\Sp_{2m_ia_i}(\kk)}(\T_i)]=\Hb_{m_i}^{a_i}$, and the claim follows. 
The claim implies \eqref{item:f-stable} by Remark \ref{rem:general-gps} \eqref{item:b}.

Let $x = sh \in \Xi$, $s\in \Sigma$, $h\in \Sigma'$. By construction 
and Lang-Steinberg's theorem
applied to $s \in \prod_i\Sp(W(m_i)^{a_i})$ and $h \in \J$, we may 
choose $g$ such that $g^{-1}F(g)=x$    as 
\begin{align}\label{eq:g=yz}
g &= yz,& &\text{where}& y&\in \prod_i\Sp(W(m_i)^{a_i}),\ 
y^{-1}F(y) = s; & z \in \J, \ z^{-1}F(z) &= h. 
\end{align}
 
Since $\Sigma'\subset C_{\G}(\Hb)\cap\J$, $\Sigma\subset C_{\G}(\J)$, 
$$(yz \J z^{-1}y^{-1})^F=(y\J y^{-1})^F=\J^F.$$
Also,  $N_i/ N_i\cap C_{\G}(u)^\circ$ is either trivial or $\simeq \Z/2$ \cite[Theorem 6.21]{LS}.
A representative of the non-trivial element is $x_i = x_i'x_i''$, where

\begin{align*}
x_i' &= \left(\begin{smallmatrix}
\id_{(a_i-1)m_i} & & & \\
&0_{m_i}&\id_{m_i}&\\
&\id_{m_i}& 0_{m_i}&\\
 & & &\id_{(a_i-1)m_i} \end{smallmatrix}\right), &
x_i'' &=
\left(\begin{smallmatrix} \id_{m_i(a_i-1)}\\
&X\\
&& \widetilde X\\
&&& \id_{m_i(a_i-1)}\end{smallmatrix}\right);
\end{align*}
here $X\in \GL_{m_i}(\kk)$ satisfies $X v X^{-1}={\Jf}_{m_i}\,^t\!v^{-1}{\Jf}_{m_i}$, and 
$\widetilde X = {\Jf}_{m_i}\,^t\!X^{-1}{\Jf}_{m_i}$. 
Then  $x_i$ normalizes each factor in  $\SL_{m_i}(\kk)^{a_i}$,  centralizes the first $a_i-1$ 
factors and  induces a non-trivial graph automorphism on the last one. 
Thus, 
\begin{align*}
(zy \Hb_{m_i}^{a_i} y^{-1} z^{-1})^F & =(z\Hb_{m_i}^{a_i} z^{-1})^F \overset{\ 
\text{Rem. \ref{rem:general-gps} \eqref{item:c}}} = z(\Hb_{m_i}^{a_i})^{Ad(x)F} z^{-1}
\\ &= z ((\Hb_{m_i}^{a_i-1})^F\times\Hb_{m_i}^{Ad(x)F})z^{-1}
\end{align*}
and the first part of \eqref{item:invaraints} follows from
Remark \ref{rem:general-gps} \eqref{item:c}. 
Finally, $N_i= N_i\cap C_{\G}(u)^\circ$ when $m_i$ is even or $m_i=1$,
hence the last restriction in \eqref{item:invaraints}. \epf

\begin{cor}\label{cor:reduction-sp}Let $u\in G$ with decomposition
\eqref{deco}. If ${\J^F}\neq1$ then  every $\Oc\in\clas{G}{u}$ 
contains a subrack that is a regular unipotent class in $\J^F$.
\end{cor}
\pf We choose the representatives of elements in $A(u)$ in $NH$ by \eqref{eq:even-symplectic-A(u)}.
For $x\in NH$, there is $g\in \Sp(\W)\times\J$ such that $g^{-1}F(g) = x$, see \eqref{eq:g=yz}. 
Then  the component $gug^{-1}_{\vert \Vc}$ 
of $gug^{-1}$ on $\Vc$ is regular in $g\J g^{-1}=\J$ and 
$\Oc_{ gug^{-1}_{\vert \Vc}}^{\J^F}$ is a subrack of $\Oc_{gug^{-1}}^{\G^F}$. \epf

\subsubsection{Preliminary results}\label{subsubsec:even-symplectic-applications}
Before starting the analysis of the various classes, we state two results needed 
for the application of Lemma \ref{lema:ACG-2-10}.

\begin{lema}\label{lem:x1x2}Let $\Oc$ be a regular unipotent class in 
either $\SL_n(q)$, $\SU_n(q)$ or $\Sp_{2n}(q)$. 
Then there are $x_1,x_2\in\Oc$ such that $(x_1x_2)^2\neq(x_2x_1)^2$.
\end{lema}
\noindent\emph{Proof.}  Let $\Jf = \Jf_n$.

\begin{step}
$\SL_n(q)$ for $n\geq2$.
By Remark \ref{rem:conjug-classes-normal-subgps}  we may assume that $\Oc \ni 
x_1=\left(\begin{smallmatrix}
1&1\\
&&\ddots&\ddots\\
&&&1&1\\
&&&&1\end{smallmatrix}\right)$; then take $x_2={\Jf}x_1{\Jf}^{-1}$. \end{step}

\begin{step}\label{case:even}
$\SU_{n}(q)$, $n$ even. By Remark \ref{rem:conjug-classes-normal-subgps} we 
assume that $\Oc\ni x_1=\left(\begin{smallmatrix}
\tu& \tv\\
0&\tu^{-1}
\end{smallmatrix}\right)$, with $\tu=\left(\begin{smallmatrix}
1&1&\dots&\dots&1\\
&&\ddots&\ddots&1\\
&&&\ddots&1\\
&&&&1\end{smallmatrix}\right)$ and $\tv = \left(\begin{smallmatrix}
1&0&\dots&0\\
%1&0&\dots&0\\
\vdots&\vdots&\vdots&\vdots\\
1&0&\dots&0\end{smallmatrix}\right)$.  Let
$x_2:=\,^t\!(\Fr_q(x_1)) = \,^t\!x_1$; we claim that $x_2\in \Oc$. 
Indeed, by definition of $\SU_{n}(q)$,  $x_2 ={\Jf}\, x_1^{-1}{\Jf} \in \Oc^{-1} = \Oc$, 
the last equality by \cite[1.4(ii)]{TZ}.
Since $x_1$ is regular, 
$C_{\SL_n(\kk)}(x_1)$ is contained in the Borel subgroup of upper triangular matrices; 
as $(x_2x_1x_2)_{21}=1$, $ x_2x_1x_2\not\in C_{\SL_n(\kk)}(x_1)$. 
\end{step}

\begin{step} $\SU_n(q)$, $n$ odd. Let $\xi\in\F_{q^2}$ satisfy $\xi^q+\xi+1=0$. 
By Remark \ref{rem:conjug-classes-normal-subgps} 
we  assume that $\Oc\ni x_1=\left(\begin{smallmatrix}
\tu&d& \xi \tv\\
0&1&b\\
0&0&\tu^{-1}
\end{smallmatrix}\right)$, with $\tu, \tv$ as in Case \ref{case:even}, $d=\left(\begin{smallmatrix}
1\\
\vdots\\
1\end{smallmatrix}\right)$,  and $b=\left(\begin{smallmatrix}
1&0&\dots&0
\end{smallmatrix}\right)$.
Now take $x_2:=\,^t\!\Fr_q(x_1)\in \Oc$ and repeat the argument for Case \ref{case:even}.
\end{step}

\begin{step} $\Sp_{2n}(q)$, $n\geq2$. There are 2 regular  classes, $\Oc$ and $\Oc'$.
Each of them is represented  by a triangular matrix whose terms 
in the upper subdiagonal are $\neq 0$. 
If $x_1$ is a representative like this, then $x_2=\sigma\trid x_1$, where 
$\sigma=\left(\begin{smallmatrix}\Jf_2&0 & 0\\
0&\id_{2n-4} & 0\\
0&0 & \Jf_2 \end{smallmatrix}\right)$, does the job.  \qed\end{step}

\begin{lema}\label{lem:y1y2}Let $n>2$ or $q>2$ and $\Oc$ a regular 
unipotent class in $\SL_n(q)$, $\SU_n(q)$, or $\Sp_{n}(q)$.
Then  there are $y_1,y_2 \in\Oc$ with $y_1\neq y_2$, $y_1y_2 = y_2y_1$.
\end{lema}
\pf By \cite[1.4(ii)]{TZ} for $\SL_n(q)$ or $\SU_n(q)$, and \cite{Gow}
for $\Sp_n(q)$, $\Oc = \Oc^{-1}$. If $n>2$ no regular element is an 
involution, so  $y_1 = y_2^{-1}$ will do.
If $n=2$ and $q>2$, then  take 
$y_1=\left(\begin{smallmatrix}1&1\\
0&1\end{smallmatrix}\right)$ and $y_2=\left(\begin{smallmatrix}1&\xi\\
0&1\end{smallmatrix}\right)$ for $1\neq \xi\in\F_q^\times$. 
\epf

\subsubsection{Analysis of the different classes}
\label{subsubsec:even-symplectic-case-by-case}
We now assume that $u\in G$ is unipotent with decomposition  
\eqref{deco}. Let $\oc$ be an arbitrary class in $\clas{G}{u}$.

\begin{lema}\label{lem:sp2-typeD} 
If \eqref{deco} contains $W(4)$, 
then $\oc$  is of type D.
  \end{lema}
\pf  By Lemma \ref{lema:representatives-a(u)-evensp}, $\oc$  
contains a subrack isomorphic to the regular  class 
in ${\SL_{4}(q)}$. Then Lemma \ref{lem:SL} applies.
\epf

\begin{lema}\label{lem:sp2-typeF}  If \eqref{deco} contains any of these terms, 
then $\oc$  is of type F.
\begin{enumerate}\renewcommand{\theenumi}{\alph{enumi}}
\renewcommand{\labelenumi}{(\theenumi)}
\item\label{item:tre}  $V(2k_j)$ with $k_j>1$ and  $q> 2$.
\item\label{item:cinque} $W(m_i)$ either with $m_i > 4$, $q > 4$; 
or else with $m_i=3$, $q > 8$; or else  with $m_i > 4$ even.
\end{enumerate}
\end{lema}
\pf \eqref{item:tre}:
 By Corollary \ref{cor:reduction-sp}, $\oc$ contains a subrack
 isomorphic to a regular class in $\Sp_{2k_j}(q)$;
 then Lemma \ref{lem:sp-regular} applies. 
\eqref{item:cinque}: By Lemma \ref{lema:representatives-a(u)-evensp},
   $\oc$ contains a subrack isomorphic to a regular class in $\SL_{m_i}(q)$, 
   or  $\SU_{m_i}(q)$ (the last occurs only when
   $m_i > 1$ is odd); 
then either Lemma \ref{lem:SL},  or else Proposition
\ref{prop:regular-char2} \eqref{item:regular-char2-psu3}, \eqref{item:regular-char2-psun}, 
or \eqref{item:regular-char2-psu4-dn}, apply.
\epf

In the next Lemma, we use that $u_{\vert W(3)}$ is regular in the image 
of $\GL_{3}(\kk)$ via \eqref{eq:image-embedding},
hence  $\Oc$ may contain a subrack isomorphic to a regular class in a 
subgroup  $\simeq \GU_{3}(8)$ when appropriate. We do so 
because the regular unipotent class in $\SU_{3}(8)$ is not known to be of type D or F.

\begin{lema}\label{lem:sp2-typeDF} If \eqref{deco} contains any of these terms, 
then $\oc$ is either of type D or else of type F. 

\begin{enumerate}\renewcommand{\theenumi}{\alph{enumi}}
\renewcommand{\labelenumi}{(\theenumi)}
\item\label{item:shest} $W(m_i)$ with $m_i>1$ odd and  $q=4$.
\item\label{item:sem} $W(3)$  and $q = 8$.
\end{enumerate}

\end{lema}

\pf We apply Lemma \ref{lema:representatives-a(u)-evensp}.
\eqref{item:shest}: $\oc$ contains a subrack isomorphic to the regular  class 
in ${\SL_{m_i}(4)}$ or in $\SU_{m_i}(4)$. Then Lemma \ref{lem:SL} 
or Proposition \ref{prop:regular-char2}
\eqref{item:regular-char2-psu3} or \eqref{item:regular-char2-psun} apply.  
\eqref{item:sem}: $\oc$ contains a subrack isomorphic to the regular class 
in ${\SL_{3}(8)}$ or in $\GU_3(8)$. Then Lemma \ref{lem:SL} or Lemma \ref{lem:GU} apply.
\epf

\begin{lema}\label{lem:sp2-product} If \eqref{deco} contains any of these terms, 
then $\oc$ is of type D.
\begin{enumerate}\renewcommand{\theenumi}{\alph{enumi}}
\renewcommand{\labelenumi}{(\theenumi)}
\item\label{item:zero}   $V(2k_i)\oplus V(2k_j)$ with $k_ik_j>1$ or $q>2$.
\item\label{item:due}  $W(m_i)\oplus V(2k_j)$ with $m_i>2$; or $m_i = 2$ 
and either $q>2$ or   $k_j>1$.
\item\label{item:uno}  $W(m_i)\oplus W(m_j)$ with either $q>2$ and 
$m_i >1$, $m_j>1$; or else $q = 2$ and  $m_i>1$ and $m_j>2$; 
or else $q = 2$,  $m_i>2$ and $m_j>1$.
\end{enumerate}
\end{lema}
\pf  By Lemma \ref{lema:representatives-a(u)-evensp}, 
$\oc$ contains a subrack $\Oc_{u_i}^{L_i} \times \Oc_{u_j}^{L_j}$ 
with each factor regular, where
$L_i\times L_j$ in each case is
\begin{align*}
&\eqref{item:zero}\hspace{-2pt}: \, \Sp_{2k_j}(q) \times \Sp_{2k_i}(q), 
\hspace{1pt} \eqref{item:due}\hspace{-2pt}: \, \SL_{m_i}(q)\times\Sp_{2k_j}(q), 
\, \text{or }\SU_{m_i}(q)\times\Sp_{2k_j}(q),
\\
&\eqref{item:uno}\hspace{-2pt}: \, \SL_{m_i}(q)\times\SL_{m_j}(q), 
\text{ or } \SL_{m_i}(q)\times\SU_{m_j}(q), \text{ or } \SU_{m_i}(q)\times\SU_{m_j}(q). 
\end{align*}
The claim follows by Lemmata \ref{lema:ACG-2-10}, \ref{lem:x1x2} and \ref{lem:y1y2}.
\epf

\begin{lema}\label{lem:sp-Wodd} If $q=2$ and 
\eqref{deco} contains a term of the  form $W(m_i)$, $m_i>1$ odd, 
then $\Oc$  is of type D.

\end{lema}
\pf {\em Step 1:} If \eqref{deco} is of the form $W(m)$ with $m>1$ 
odd, then $\Oc$  is of type D.
Indeed, there are two classes of this type \cite[Theorem 6.21]{LS}. 
One class contains a subrack isomorphic to the regular unipotent 
class in $\GU_m(2)$, and we invoke Lemma \ref{lem:GU}.
We consider next the second class. Assume first that $m=3$.
Then this class is represented by
$v = \left(\begin{smallmatrix}
1&1&1\\
&1&1\\
&&1\\
&&&1&1&0\\
&&&&1&1\\
&&&&&1\end{smallmatrix}\right)=x_{\alpha_1}(1)x_{\alpha_2}(1)$.
\begin{align*}
\text{Let }&s_2:= \left(\begin{smallmatrix}
1\\
&0&1\\
&1&0\\
&&&0&1\\
&&&1&0&\\
&&&&&1\end{smallmatrix}\right), &r=s_2\trid v=\left(\begin{smallmatrix}
1&1&1\\
&1&0\\
&1&1\\
&&&1&0&1\\
&&&1&1&0\\
&&&&&1\end{smallmatrix}\right)=x_{\alpha_1+\alpha_2}(1)x_{-\alpha_2}(1),\\
&s_3=\left(\begin{smallmatrix}
\id_2\\
&0&1\\
&1&0\\
&&&\id_2\end{smallmatrix}\right),
& s=s_3\trid v = \left(\begin{smallmatrix}
1&1&0&1\\
&1&0&1\\
&&1&0&1\\
&&&1&0&0\\
&&&&1&1\\
&&&&&1\end{smallmatrix}\right)=x_{\alpha_1}(1)x_{\alpha_2+\alpha_3}(1).
\end{align*}
Then $(rs)^2\neq(sr)^2$ by \eqref{eq:chev}. In addition, $r,s\in\Pa^F$ 
where $\Pa$ is the standard parabolic subgroup of $\G$ associated 
with the simple root $\alpha_2$ so
Remark \ref{rem:parabolic} (1) applies.
Assume $m>3$. Then the class is represented by
\begin{align*}v=\left(\begin{smallmatrix}
1&1&\cdots&1\\
&1&\cdots&1\\
&&1&1\\
&&&1\\
&&&&1&1&0&\cdots&0\\
&&&&&1&1&0&\cdots\\
&&&&&&1&1&0\\
&&&&&&&1&1\\
&&&&&&&&1\end{smallmatrix}\right)=x_{\alpha_1}(1)x_{\alpha_2}(1)\cdots x_{\alpha_{n-1}}(1).\end{align*}
We apply Lemma \ref{lem:parabolic} to $\Pa^F$, where  $\Pa$ is the 
standard parabolic associated with the simple 
roots $\alpha_n,\alpha_{n-1},\alpha_{n-2}$, using the case $m=n=3$.

{\em Step 2:} We now prove the Lemma. 
Let $u_i$ be the component of $u$ in $\M_i:=\Sp(W(m_i))$. Choosing 
each representative $x$ of $A(u)$ in $\Xi$ and 
the corresponding element $g$ as in Lemma \ref{lema:representatives-a(u)-evensp},  
we have  $g\M_ig^{-1}=\M_i$,
 and $(g\M_i g^{-1})^F=g \M_i^{Ad(x)\circ F}g^{-1}$. So,
$\Oc_{gug^{-1}}^G$ will contain a subrack isomorphic to
$\Oc_{u_i}^{\M_i^{Ad(x)\circ F}}$.  The component in $\M_i$ of 
each term in $\Xi$  is either trivial or, possibly, $x_i'x_i''$, with notation as in Lemma 4.6.
Then, the two possible subracks are isomorphic to those in Case 1, whence the statement.
\epf

\begin{obs}\label{remark:even-sp-1st-reduction}
By the previous Lemmata, it remains to consider the following forms of \eqref{deco}, see  
\ref{tab:unipotent-chevalley-collapse-symplectic-even-I} for details:
\begin{gather}\label{deco>2}
\begin{tabular}{c c c}
$q>2$: & $V =  W(1)^{a}\oplus W(2)$,  & $0 \le a$
\\  & $V = W(1)^{a} \oplus V(2)$, & $0 \le a$
\\
$q =2$: &  $V =  W(1)^{a}\oplus W(2)^{b} \oplus V(2)^{c}$,
 &   $0 \le a,b$; $0\le c \le 2$
 \\ & $V = W(1)^{a} \oplus V(2k)$ & $0 \le a$; $1 < k$
 \end{tabular}
\end{gather}
\end{obs}

\begin{gather} \label{tab:unipotent-chevalley-collapse-symplectic-even-I} 
\stepcounter{tabla}\tag*{Table \thetabla}
\begin{tabular}{|c|c|c|c|c|}
\hline  \eqref{deco} $\supseteq$ & $k_j$ & $m_i$ & $q$ &  Criterium  \\
\hline  \hline
 $V(2k_j)$ & $>1$ & -- & $>2$   & F,  \ref{lem:sp2-typeF} \eqref{item:tre} \\
 \hline
 $W(m_i)$ & -- & $>4$ & $>4$   & F,  \ref{lem:sp2-typeF} \eqref{item:cinque} \\
 \cline{3-5}
 & & $>4$ even & all  & F,  \ref{lem:sp2-typeF} \eqref{item:cinque}  \\
  \cline{3-5}
 &  & $> 4$ odd & 4   &  \small{F or D,}   \ref{lem:sp2-typeDF}  \\
  \cline{3-5}
 &  & 4 & all & D,  \ref{lem:sp2-typeD}\\
 \cline{3-5}
 &  & 3 &  $>8$     & F,  \ref{lem:sp2-typeF} \eqref{item:cinque} \\
\cline{4-5}
& & & 8, 4   &  \small{F or D,}  \ref{lem:sp2-typeDF} \\
  \cline{3-5}
 &  & $> 1$ odd & 2   &  D,  \ref{lem:sp-Wodd} \\
  \hline 
\small{$V(2k_i)\oplus V(2k_j)$} &  &  & $>2$& D, \ref{lem:sp2-product} \eqref{item:zero}
\\ \cline{2-4}
 & \small{$k_ik_j>1$} &  & $2$& 
\\  \hline 
\small{$W(m_i)\oplus V(2k_j)$} & & $>2$ & all &
\\ \cline{2-4}
 & &  $2$ & $>2$ &  D, \ref{lem:sp2-product} \eqref{item:due} 
 \\ \cline{2-4}
 & $>1$ &  $2$ & $2$ & 
 \\  \hline 
\small{$W(m_i)\oplus W(m_j)$} & & \small{$m_i ,m_j>1$} & $> 2$ &
\\ \cline{2-4}
 & &  \small{$m_i>1$, $m_j>2$} & $2$ &  D, \ref{lem:sp2-product} \eqref{item:uno} 
 \\ \cline{2-4}
 &  &  \small{$m_i>2$, $m_j>1$} & $2$ & 
 \\
 \hline 
 \end{tabular}
  \end{gather}

Recall the conventions in \eqref{deco} on $a_i$, $b_j$.

\begin{lema}\label{lem:sp2-a>0} If  \eqref{deco} is of either of the following forms,
then $\clas{G}{u}$ consists of only one class which is of type  D:
\begin{enumerate}\renewcommand{\theenumi}{\alph{enumi}}\renewcommand{\labelenumi}{(\theenumi)}
\item\label{item:undici}  $W(1)^{a_1}\oplus W(2)^{a_2}$ or $W(2)^{a_2}$,  $a_2 > 1$.
\item\label{item:duodici}  $W(1)^a \oplus W(2)^{a_2} \oplus V(2)^{b_1}$, $0\le a$.
\end{enumerate}
\end{lema}
\pf In all cases $\clas{G}{u}$ has only one class $\oc$ by \cite[Theorem 6.21]{LS}. 

\eqref{item:undici}:
Let $x=\left(\begin{smallmatrix}
        1&1\\
        0&1\end{smallmatrix}\right)$.
Let $v$ be the block diagonal matrix
        $$(x\otimes\id_{a_2},\id_{2a_1},\Jf_{2}\,^tx^{-1}\Jf_{2}\otimes\id_{a_2});$$
        then $v\in \Oc$ because its decomposition \eqref{deco} is $W(1)^{a_1}\oplus W(2)^{a_2}$. 
       Now $v$ lies in the subgroup $H\simeq \SL_n(q)$ of 
        matrices $\left(\begin{smallmatrix}
y & 0\\ 0& \Jf_{n} \,^{}ty^{-1}\Jf_{n}\end{smallmatrix}\right)$, with $y\in\SL_n(q)$. 
If $a_2>1$, then
 $\Oc$ contains a subrack isomorphic to a unipotent class of type $(2,...,2)$ ($a_2$ times) in $\SL_{2a_2}(q)$, 
and Lemma \ref{lem:SL} applies.
        
%Assume $a_2=1$. If $a_1>2$, then $v$ lies in $H$ and $\Oc$ contains a subrack isomorphic to a unipotent class of 
%type $(2,1,1,1,)$ in $\SL_n(q)$, and Lemma \ref{lem:SL} applies.

\smallbreak
\eqref{item:duodici}: It is enough to consider $W(2) \oplus V(2)$. Let 
$v=\left(\begin{smallmatrix}
1&0&0&0&0&1\\
&1&1&0&0&0\\
&&1&0&0&0\\
&&&1&1&0\\
&&&&1&0\\
&&&&&1\end{smallmatrix}\right)=x_{\alpha_2}(1)x_{2\alpha_1+2\alpha_2+\alpha_3}(1) \in \Oc.$
We set \begin{align*}
\sigma&= \left(\begin{smallmatrix}
0&1&0\\
0&0&1\\
1&0&0\\
&&&0&0&1\\
&&&1&0&0\\
&&&0&1&0\end{smallmatrix}\right),& 
&z = \left(\begin{smallmatrix}
1&0&0&0&0&0\\
&1&0&1&1&0\\
&&1&0&1&0\\
&&&1&0&0\\
&&&&1&0\\
&&&&&1\end{smallmatrix}\right),&  
y &= \left(\begin{smallmatrix}
1&0&0&0&0&0\\
&1&0&0&0&0\\
&&1&1&0&0\\
&&&1&0&0\\
&&&&1&0\\
&&&&&1\end{smallmatrix}\right).
\end{align*}
and 
\begin{align*}
&r=(z\sigma)\trid v=\left(\begin{smallmatrix}
1&1&0&1&1&0\\
0&1&0&0&0&1\\
0&0&1&1&0&1\\
0&0&0&1&0&0\\
0&0&0&0&1&1\\
0&0&0&0&0&1
\end{smallmatrix}\right),  &s=y\trid v=\left(\begin{smallmatrix}
1&0&0&0&0&1\\
0&1&1&1&0&0\\
0&0&1&0&1&0\\
0&0&0&1&1&0\\
0&0&0&0&1&0\\
0&0&0&0&0&1
\end{smallmatrix}\right)
\end{align*} 
The discussion in \cite[3.1]{ACG-I} implies that 
$\Oc_r^{\U^F}\neq \Oc_s^{\U^F}$. A direct computation 
shows that $(rs)^2\neq(sr)^2$.\epf 

\begin{lema}\label{lem:sp2-a>1} If  \eqref{deco} is equal to 
$W(1)^{a_1}\oplus W(2)$,  $a_1 > 1$
then $\clas{G}{u}$ consists of only one class which is of type  F.
\end{lema}
\pf In all cases $\clas{G}{u}$ has only one class $\oc$ by \cite[Theorem 6.21]{LS}. 
It is enough to prove the statement for $a_1=2$. 
Let $x$  and $v$ be as in Lemma \ref{lem:sp2-a>0} (a). 
Then $v\in\Oc$ because it has decomposition \eqref{deco} equal to $W(2)\oplus W(1)^2$.
We consider the following elements of  $G$ 
\begin{align*}
\sigma&= \left(\begin{smallmatrix}
0&0&1&0\\
1&0&0&0\\
0&1&0&0\\
0&0&0&1\\
&&&&1&0&0&0\\
&&&&0&0&1&0\\
&&&&0&0&0&1\\
&&&&0&1&0&0\end{smallmatrix}\right),& 
&\tau =  \left(\begin{smallmatrix}
1&0&0&0\\
0&0&0&1\\
0&1&0&0\\
0&0&1&0\\
&&&&0&1&0&0\\
&&&&0&0&1&0\\
&&&&1&0&0&0\\
&&&&0&0&0&1\end{smallmatrix}\right),&  
\omega&=\left(\begin{smallmatrix}
\id_2\\
&0&\id_2\\
&\id_2&0\\
&&&\id_2
\end{smallmatrix}\right)
\end{align*}
and the following elements of $\Oc$:
\begin{align*}
&r_1=v, &r_2=\sigma\trid r_1=(1,x,\id_2,x,1)\\
&r=\tau\trid r_2=(\id_2,x\otimes \id_2,\id_2),& r_3=(r_2\tau)\trid r_2 =\left(\begin{smallmatrix}
1&0&0&0\\
0&1&0&1\\
0&0&1&1\\
0&0&0&1\\
&&&&1&1&1&0\\
&&&&0&1&0&0\\
&&&&0&0&1&0\\
&&&&0&0&0&1\\
\end{smallmatrix}\right)\\
&r_4=(r_1r\omega)\trid r_2
=\left(\begin{smallmatrix}
1&0&0&0&1&1&0&0\\
0&1&0&0&1&1&0&0\\
0&0&1&0&0&0&1&1\\
0&0&0&1&0&0&1&1\\
&&&&1&0&0&0\\
&&&&0&1&0&0\\
&&&&0&0&1&0\\
&&&&0&0&0&1\\
\end{smallmatrix}\right).
\end{align*} 

Then, $H:=\langle r_1, r_2,r_3, r_4\rangle\subset \U^F$ and  $\Oc_{r_i}^{\U^F}\neq\Oc_{r_j}^{\U^F}$, for $i\neq j$, hence  $\Oc_{r_i}^{H}\neq\Oc_{r_j}^{H}$. A direct computation shows that $r_i\trid r_j\neq r_j$ for $i\neq j$.
\epf

\begin{lema}\label{lem:missingclass}Assume $q>2$. 
If  \eqref{deco} is equal to $W(1)\oplus W(2)$,  
then $\clas{G}{u}$ consists of only one class which is of type  F.
\end{lema}
\pf There is only one class $\Oc$ with \eqref{deco} equal to $W(2)\oplus W(1)$, which is represented by
$r_1=\id_{6}+(e_{2,3}+e_{4,5})$, \cite[Theorem 6.21]{LS}.

Let $\zeta\in F_q^\times\setminus 1$ and let us consider the following elements of $G$:
\begin{align*}
s_1:=\left(\begin{smallmatrix}
0&1&0\\
1&0&0\\
0&0&1\\
&&&1&0&0\\
&&&0&0&1\\
&&&0&1&0
\end{smallmatrix}\right)&& s_2:=\left(\begin{smallmatrix}
1&0&0\\
0&0&1\\
0&1&0\\
&&&0&1&0\\
&&&1&0&0\\
&&&0&0&1
\end{smallmatrix}\right)&& s_3:=\left(\begin{smallmatrix}
\id_2\\
&0&1\\
&1&0\\
&&&\id_2
\end{smallmatrix}\right)
\end{align*}

We construct the following elements in  $\Oc$:
\begin{align*}
&r_2=(r_1 s_2 s_1)\trid r_1=\id_{6}+(e_{1,2}+e_{5,6})+(e_{1,3}+e_{4,6})\\
&r_3=\left((\id_6+(e_{2,1}+e_{6,5}))s_3s_1\right)\trid r_1=\id_{6}+(e_{1,4}+e_{3,6})+(e_{2,4}+e_{3,5})\\
&r_4=\left({\rm diag}(1,\zeta,1,1,\zeta^{-1},1)(\id_{6}+e_{3,4})(\id_{6}+e_{1,2}+e_{5,6})\right)\trid r_1=
\left(\begin{smallmatrix}
1&0&1&1&0&0\\
&1&\zeta&\zeta&0&0\\
&&1&0&\zeta&1\\
&&&1&\zeta&1\\
&&&&1&0\\
&&&&&1
\end{smallmatrix}\right)
\end{align*}
A direct computation shows that $r_i\trid r_j\neq r_j$ for $i\neq j$. Moreover, as $H:=\langle r_1,r_2,r_3,r_4\rangle\subset \U^F\subset \SL_{6}(q)$, the usual argument shows that
$\Oc_{r_i}^H\neq \Oc_{r_j}^H$ for $i\neq j$. \epf

\begin{lema}\label{lem:k2q2}
If $q=2$ and \eqref{deco} is of the form   $W(1)^a\oplus V(4)$, $0\le a$, 
then $\oc$ is of type D.
\end{lema}

\pf There are $2$ classes like this  \cite[Theorems 6.6, 6.12]{LS}; both 
contain a subrack isomorphic to one of the regular  classes in 
 $\Sp_4(2)\simeq\Sim_6$, which  corresponds either to the partition $(4,2)$ 
 or else to $(4,1^2)$. These  are  of type D by 
 \cite[4.1]{AFGV-ampa}.
\epf

\begin{lema}\label{lem:sp-regular-3}
If  $q=2$ and  \eqref{deco} contains the  term $V(2k)$, $k\geq 3$,
then $\oc$ is of type  D.

\end{lema}

\pf By Corollary \ref{cor:reduction-sp},  $\oc$ contains a subrack isomorphic to
one of the two regular unipotent classes in $\Sp_{2k}(2)$. So, we 
may assume $k = n$ for notational purposes.
For both classes, there is a representative $v$ lying in 
$x_{\alpha_1}(1)\cdots x_{\alpha_n}(1)\U'$, $\U'$ as in 
Proposition \ref{prop:regular}. 
Let $\Pa$ be the standard parabolic subgroup of $\G$ corresponding to the 
simple roots $\alpha_{n-1}$ and $\alpha_n$, and let $\Le$ be  its 
standard Levi subgroup, whose derived subgroup is isomorphic to $\Sp_4(\kk)$. 
Recall the notation in Remark \ref{rem:parabolic}.
Then  $v\in\U^F< \Pa^F$ and if $v = v_L v_P$ is its decomposition 
according to $\Pa^F=\Le^F\ltimes\urad^F$, then $u_L$ is regular unipotent in $\Le^F$. 
The result follows from Lemmata \ref{lem:parabolic} and \ref{lem:k2q2}.
\epf

\begin{lema}\label{lem:22q2}If $q=2$ and \eqref{deco} is of the form  
$W(1)^{a_1}\oplus V(2)^2$, then $\Oc$  is of type D.
\end{lema}
\pf There is only one class in $\clas{G}{u}$ by \cite[Theorem 6.21]{LS}.
We may assume $a_1 = 1$, $n=3$ as in the proof of Lemma \ref{lem:sp2-a>0} \eqref{item:undici}. 
Then $\Oc$ is represented by $xy=x_{2(\alpha_1+\alpha_2)+\alpha_{3}}(1)x_{2\alpha_2+\alpha_3}(1)$. 
It contains the subrack
$X\times Y$ for $X=\Oc_{x}^H$, $Y=\Oc_{y}^K$ with $H\simeq \SL_2(2)$ 
being the subgroup corresponding to the root 
$2(\alpha_1+\alpha_2)+\alpha_{3}$ and $K\simeq \Sp_4(2)\simeq\Sim_6$ 
the subgroup corresponding to the roots $\alpha_2$ and $\alpha_3$. 
Since all conjugacy classes of involutions in $\Sim_6$ contain distinct commuting 
elements, we apply Lemmata \ref{lema:ACG-2-10} and \ref{lem:x1x2}.
\epf

\subsection{Proof of Theorem \ref{th:unipotent-chevalley-collapse}}\label{subsec:th1-symplectic}
We first study classes that do not collapse. 

\begin{lema}\label{lem:sp-hook} Let $u\in G$ unipotent with partition $(1^{2n-2},2)$. 

\begin{enumerate}
 \item\label{item:hook-odd} If $q$ is odd and either not a square or $9$, 
 then $\clas{G}{u}$ consists of two cthulhu classes.

  \item\label{item:hook-even} If $q$ is even, then $\clas{G}{u}$ 
  consists of a unique cthulhu class.
 \end{enumerate}
 \end{lema}
 
\pf   If $q$ is even, the decomposition \eqref{deco} of any element with Jordan form $(1^{2n-2},2)$,  is 
necessarily $W(1)^{2n-2}\oplus V(2)$. Thus we have only one conjugacy class in $\G$ with this form.
  
For any $q$, we fix $u =\left(\begin{smallmatrix}1&0& 1\cr
0&\id_{2n-2}&0\\
0&0&1\end{smallmatrix}\right) =x_{\beta}(1)$, where $\beta\in \Phi^+$ is the highest root. 
If $q$ is even, $\clas{G}{u}$ has a unique class by \cite[Theorems 6.6, 6.12]{LS}. 
If $q$ is odd, then, by \cite[Theorem 3.1(v)]{LS} and arguing as in 
Lemma~\ref{lem:sp-2columns}, $\clas{G}{u}$  consists of two classes 
represented by $u$ and $\left(\begin{smallmatrix}1&0&\zeta\cr
0&\id_{2n-2}&0\\
0&0&1\end{smallmatrix}\right) =x_{\beta}(\zeta)$, for $\zeta\in \F_q^\times$ not a square. 
We show (for any $q$) that every subrack of $\Oc_u^{G}$ generated by two elements
is either abelian or indecomposable, implying that $\Oc_u^{G}$ is cthulhu. 
The same argument applies to the other class, when $q$ is odd.

Assume there is $g\in G$ such that $v=gug^{-1}\in \Oc_u^{G}$ 
and $uv\neq vu$. We claim that 
the rack generated by $u$ and $v$ 
is indecomposable.
Consider the Bruhat decomposition $g=y tn_w z$, with $y,z\in \U^F$, $t\in \T^F$ 
and $n_w\in N_G(\T)$ with class $w\in W$. 
By  \eqref{eq:chev}, $u\in Z(\U^F)$, so that $v=huh^{-1}$ with $h=y tn_w$. 
Now the subrack generated by $u$ and $v$ is isomorphic
to the the subrack generated by $u = y^{-1}uy$ and $y^{-1}vy = kuk^{-1}$ 
with $k=tn_w$, so we may assume that $v = kuk^{-1}$.
Now a direct computation gives $v = x_{w\beta}(\eta)$ for some  
$\eta\in\F_q^\times$ \cite[Theorems 24.10; 8.17(e)]{MT}.
The assumption $uv\neq vu$ forces $w\beta\in -\Phi^+$ and $w\beta+\beta\in\Phi\cup\{0\}$. 
As the root system is of type $C_n$, this is possible only if $w\beta=-\beta$. 
An  element in $N_G(\T)$ mapping $u= x_{\beta}(1)$ to $v= x_{-\beta}(\eta)$ is of the form
$$\left(\begin{smallmatrix}
0&0&\xi\\
0&X&0\\
-\xi^{-1}&0&0
\end{smallmatrix}\right)=
\left(\begin{smallmatrix}
\xi&0&0\\
0&\id _{2n-2}&0\\
0&0&\xi^{-1}
\end{smallmatrix}\right)
\left(\begin{smallmatrix}
0&0&1\\
0&\id_{2n-2}&0\\
-1&0&0
\end{smallmatrix}\right)\left(\begin{smallmatrix}
1&0&0\\
0&X&0\\
0&0&1\end{smallmatrix}\right)=\beta^\vee(\xi)n_{\beta}Y,
$$ for $X\in\Sp_{2n-2}(q)$ and $\xi\in \F_q^\times$. Hence $\eta=-\xi^{-2}$ 
and $Y$ commutes with $u$.

Let $H = \left\{\ \left(\begin{smallmatrix}
a&&b\\
&\id_{2n-2}\\
c&&d\end{smallmatrix}\right) \in G: ad-bc=1 \right\}\simeq\SL_2(q)$.
Then $u, v, \beta^\vee(\xi)$, $n_\beta\in H$ and $v\in \Oc_u^H$. 
By \cite[Lemma 3.5]{ACG-I}, $\Oc_u^H$ is sober, hence the rack generated 
by $u$ and $v$ is indecomposable. \epf

\begin{obs}\label{obs:cthulhu-symplectic}
Let $q$ be odd, $u=x_{\beta}(1)$, $\syp_n = \Oc_u^G$ and let $u'= x_{\beta}(\zeta)$, 
for $\zeta\in \F_q^\times$ not a square. Then
$\syp_n \simeq \Oc_{u'}^G$ as racks because the outer automorphism of 
$G = \Sp_{2n}(q)$ given by conjugation
by the matrix $\diag(\id_n, \zeta^{-1}\id_n)$ maps $u$ to $u'$.
Thus for $q$ even or $q=9$, or $q$ odd and not a square, we have a family of 
cthulhu racks $(\syp_n)_{n\in \N}$, 
with $\syp_1$ the sober rack $\oc^{\SL_2(q)}_x$ with $x$ nontrivial unipotent.
Note that $\syp_n \subset \syp_{n+1}$ and 

\begin{align}
\vert \syp_n\vert  = \begin{cases}\frac{(q^{2n}-1)}{2}, &\text{if $q$ is odd,}\\
 (q^{2n}-1), &\text{if $q$ is even.}\\
                     \end{cases}
\end{align}
\end{obs}

\begin{lema}\label{lem:22q2-1}If $q=2$ and \eqref{deco} is of the form  
$V(2)^2$, then $\Oc$  is cthulhu.
\end{lema}

\pf There is only one class in $\clas{G}{u}$ by \cite[Theorem 6.21]{LS}.
Here $\G^F=\Sp_4(2)\simeq\Sim_6$ and $\oc$ corresponds to the partition $(1^2,2^2)$.
By \cite[Remark 4.2 (e)]{AFGV-ampa}, $\oc$ is not of type D.  We will show 
that it cannot be of type F either. 
For $i\in\I_4$, let $r_i\in\Oc$ with $[r_i,r_j]\neq1$ and 
$\Oc_{r_i}^{\langle r_i,r_j\rangle}\neq\Oc_{r_j}^{\langle r_i,r_j\rangle}$ for $i\neq j$. 
Then  for every $i\neq j$, the permutations $r_i$ and $r_j$ may not have a 
$2$-cycle in common, and $\langle r_i,r_j\rangle$ cannot
be contained in a standard subgroup isomorphic to $\Sim_4, \Sim_5,$ or 
$\Sim_3\times \Sim_3$. 
If $r_4=(12)(34)$, then for $i\in \I_3$ we necessarily have $r_i$ either in 
$A=\{(13)(56), (14)(56), (23)(56), (24)(56)\}$ or in 
$B=\{(15)(26), (16)(25), (35)(46), (45)(36)\}.$ However, if $r_2 \in A$,
respectively $B$, then $r_3, r_4$ must lie in $B$, respectively 
$A$, leading to a contradiction.
\epf

\begin{lema}\label{lem:W(2)W(1)not-D}Assume $q$ is even.
If  \eqref{deco} is equal to 
$W(1)^{a_1}\oplus W(2)$,  
then $\clas{G}{u}$ consists of only one class $\Oc$ which is not of type  D. If $a_1=1$ and  $q=2$, then $\Oc$ is cthulhu.
\end{lema}
\pf  $\clas{G}{u}=\{\Oc\}$ by \cite[Theorem 6.21]{LS}. We shall prove that for any two elements $r,\,s \in\Oc$ such that $(rs)^2\neq(sr)^2$, it holds $\Oc_r^{\langle r,s\rangle}=\Oc_s^{\langle r,s\rangle}$. Let $\gamma=\varepsilon_1+\varepsilon_2$ be the highest short root in the root system of $\G$. The class $\Oc$ is represented by $r=x_{\gamma}(1)=\id_{2n}+e_{1,2n-1}+e_{2,2n}$, which is central in $\U^F$ by \eqref{eq:chev} and \ref{tab:c11-alpha-beta}. Let $s=g\trid r\in \Oc$ satisfy $(sr)^2\neq (rs)^2$ and
let  $g=u\dw v\in\U^F N_G(\T) \U^F$ be the Bruhat decomposition of $g$. Then
$s=(u \dw)\trid r=u\trid x_{w(\gamma)}(\eta)$ for some $\eta\in\F_q^\times$. Conjugating by $u^{-1}$ we may assume $s=x_{w(\gamma)}(\eta)$. Now, as $sr\neq rs$, we necessarily have
$w(\gamma)\in\{-\gamma, -\varepsilon_1\pm\varepsilon_k,-\varepsilon_2\pm\varepsilon_k, k\neq1,2\}$. We claim that $w(\gamma)=-\gamma$. Assume indeed that 
$w(\gamma)\in  \{-\varepsilon_1\pm\varepsilon_k,-\varepsilon_2\pm\varepsilon_k, k\neq1,2\}$. By \eqref{eq:chev}, we have $rsrs\in \U_{\gamma+w(\gamma)}$, so it is an involution, leading to a contradiction. Thus, 
$$H:=\langle r,s\rangle\simeq \langle\left(\begin{smallmatrix}
1&1\\
0&1\end{smallmatrix}\right),\,\left(\begin{smallmatrix}
1&0\\
\eta&1\end{smallmatrix}\right)\le\SL_2(q).$$
Since the non-trivial unipotent rack in $\SL_2(q)$ is sober, we have the first statement. The second one follows from a computation with GAP.
\epf

\begin{lema}\label{lem:sp2-(010)}Let $G=\Sp_{4}(q)$ for $q$ even and let $\Oc_u^{G}$ 
be a class corresponding to $W(2)$.
Then $\clas{G}{u}$ contains  a unique class which is cthulhu.
\end{lema}
\pf  The root system of $\G = \Sp_{4}(\kk)$ is of type $C_2$, so there exists a non-standard 
graph automorphism $\theta$ 
interchanging long and short roots \cite[12.1]{carter-simple}, commuting with $F$. 
Thus, $\theta$ induces an automorphism on $\G^F$
mapping the class of type $W(1)^2\oplus V(2)$, represented by $x_{\alpha_1}(1)$, 
onto the class of type $W(2)$, represented by 
$x_{\alpha_2}(1)$. The claim follows from  Lemma  \ref{lem:sp-hook} \eqref{item:hook-even}. \epf

We next show that the classes not listed in \ref{tab:uno} collapse. 
Let $\oc$ be a unipotent class in $G$. 
We summarize in \ref{tab:unipotent-chevalley-collapse-symplectic-odd} 
the results in \S \ref{subsec:odd-symplectic} 
proving the claim for $q$ odd. 

\begin{gather} \label{tab:unipotent-chevalley-collapse-symplectic-odd} 
\stepcounter{tabla}\tag*{Table \thetabla}
\begin{tabular}{|c|c|c|}
\hline  $q$, $n$ & type  $(1^{r_1}, 2^{r_2}, \dots, n^{r_n})$ & Criterium  \\
\hline  \hline
& $\exists\, i >3: r_i\neq 0$   & type D,  \ref{lem:sympl-odd}   \\ \hline
  $> 9$   square    & $(1^{r_1}, 2^{r_2})$, $r_2 > 0$ & type D,   \ref{lem:sympl-odd} \\
  \hline
 $q>3$, or $n>2$ & $(1^{r_1}, 2^{r_2})$, $r_2 > 1$ & type D,  \ref{lem:sp-2columns}
\\  \hline
   & $(1^{r_1}, 3^{r_3})$, $r_3 > 0$ &  type D,  \ref{lem:sp31} \\ \hline
   & $(1^{r_1}, 2^{r_2}, 3^{r_3})$, $r_2r_3 > 0$ & type D,  \ref{lem:321sp}  \\\hline 
  3 & $(2^2)$ &  one  of type D,  \ref{lem:sp43}\\
  \hline
 \end{tabular}
\end{gather}

\bigbreak
Assume that $q$ is even. We show how the results in  \S
\ref{subsec:even-symplectic} imply the claim. By Remark \ref{remark:even-sp-1st-reduction} 
we may assume
that \eqref{deco} has the form \eqref{deco>2}. For all $q$ even, we have
\begin{itemize}[leftmargin=*]\renewcommand{\labelitemi}{$\diamond$}
\item $V =  W(2)$: cthulhu,  Lemma \ref{lem:sp2-(010)}.
  \item $V =  W(1)^{a}\oplus V(2)$,  $0 \le a$: cthulhu,  Lemma \ref{lem:sp-hook} (2).
\end{itemize}

\medbreak
\emph{Case $2 < q$ even.} The remaining cases are disposed as follows.
\begin{itemize}[leftmargin=*]\renewcommand{\labelitemi}{$\diamond$}
 \item $V =  W(1)^{a}\oplus W(2)$,  $1 \leq a$: type F,  
 Lemmata \ref{lem:sp2-a>1} and \ref{lem:missingclass}.  
\end{itemize}

\medbreak
\emph{Case $q=2$}. Here we invoke the following statements.
\begin{itemize}[leftmargin=*]\renewcommand{\labelitemi}{$\diamond$}
 \item $V =   W(1)^{a}\oplus W(2)^{b} \oplus V(2)^{c}$,  $0 < bc$: type D,  
 Lemma \ref{lem:sp2-a>0} \eqref{item:duodici}.
  \item $V =   W(1)^{a}\oplus W(2)^{b}$, $1 <ab$; or $W(2)^{b}$, $1 < b$: 
  type D or F,  Lemmata \ref{lem:sp2-a>0} \eqref{item:undici}, \ref{lem:sp2-a>1}.
  \item $V =  W(1)^{a}\oplus V(2)^2$,  $0 < a$: type D,  Lemma \ref{lem:22q2}.
    \item $V =  V(2)^2$: cthulhu,  Lemma \ref{lem:22q2-1}.
  \item $V = W(1)^{a} \oplus V(2k)$, $0 \le a$, $k \ge 2$: type D, 
  Lemmata \ref{lem:k2q2}, \ref{lem:sp-regular-3}. 
  \item  $V= W(2)\oplus W(1)$, cthulhu,  Lemma \ref{lem:W(2)W(1)not-D}.
\end{itemize}


\begin{thebibliography}{AFGV2}%{9991}

\bibitem[ACG]{ACG-I}  N. Andruskiewitsch, G. Carnovale, G. A.
Garc\'ia.
\emph{Finite-dimensional pointed Hopf algebras over finite simple groups of Lie type  I. 
Non-semisimple classes in $\PSL_n(q)$}, J. Algebra, to appear.


\bibitem[AFGV1]{AFGV-ampa} N. Andruskiewitsch, F. Fantino, M.
Gra\~na and  L. Vendramin, 
\emph{Finite-di\-mensional pointed
Hopf algebras with alternating groups are trivial},
Ann. Mat. Pura Appl. (4),  \textbf{190}  (2011), 225--245.

\bibitem[AFGV2]{AFGV-espo} \bysame,
\emph{Pointed Hopf algebras over the sporadic simple groups}.
J. Algebra \textbf{325} (2011), pp. 305--320.


\bibitem[AG]{AG-adv} N. Andruskiewitsch and M. Gra\~na,
    \emph{From racks to pointed Hopf algebras},
Adv. Math.  \textbf{178}  (2003), 177--243.

\bibitem[Bou]{Bou} {N. Bourbaki},
{\it Groupes et alg\`{e}bres de Lie Chap. IV,V,VI}, Hermann, Paris, 1968.

%\bibitem[Bo]{borel} A. Borel,
%{\it Linear algebraic groups}, Second edition. Graduate Texts in Mathematics, 126. Springer-Verlag, New York, 1991.

\bibitem[C]{carter-simple} R. W. Carter, 
\textit{Simple groups of Lie type}.
John Wiley \& Sons, Ltd., 1972. 

%\bibitem[C2]{carter-centralizers} \bysame \emph{Centralizers of semisimple elements in finite groups of Lie type}.
%Proc. London Math. Soc. (3) \textbf{37} (1978), 491–507.


%\bibitem[C3]{carter-classical} \bysame,
%\emph{Centralizers of semisimple elements in the finite classical groups}. Proc. London Math. Soc. (3) 42 (1981), no. 1, 1–41.


%\bibitem[C4]{carter} \bysame,
%\textit{Finite groups of Lie type. Conjugacy classes and complex characters}.
%Reprint of the 1985 original. John Wiley \& Sons, Ltd., Chichester, 1993.

%\bibitem[Ch]{chang}B. Chang,
%\emph{The conjugate classes of Chevalley groups of type $G_2$,}
%J. Algebra \textbf{9}  (1968), 190-211.

%\bibitem[CH]{CH} M. Cuntz and I. Heckenberger,
%\emph{Finite Weyl groupoids}. \texttt{arXiv:1008.5291v1}.


%\bibitem[E]{eno}{H. Enomoto},
%\emph{The conjugacy classes of Chevalley groups of type $G_2$ over finite fields of characteristic $2$ or $3$,}
%J. Fac. Sci. Univ. Tokyo, Sect. I,  \textbf{16}  (1970), 497-512.

%\bibitem[FaV]{FV} {F. Fantino} and {L. Vendramin},
%\textit{On twisted conjugacy classes of type D in sporadic simple groups}.
% Contemp. Math. \textbf{585} (2013) 247-259.


\bibitem[Go]{Gow}{R. Gow},
\emph{Products of two involutions in classical groups of characteristic 2,}
J. Algebra
\textbf{71}, 583-591 (1981).


%\bibitem[HV]{HV} I. Heckenberger and L. Vendramin,
%\emph{Nichols algebras over groups with finite root system of rank two II}.  \texttt{http://arxiv.org/abs/1302.0213}.



\bibitem[Hu]{Hu} J. E. Humphreys,
\emph{Conjugacy classes in semisimple algebraic groups},
Amer. Math. Soc., Providence, RI, 1995.


%\bibitem[LaT]{LaT}R. Lawther, D. Testerman,
%\emph{Centres of Centralizers of Unipotent Elements in Simple Algebraic Groups,}
%Mem. Amer. Math. Soc. Vol. 210, Number 988, (2011).

\bibitem[LS]{LS}M. W. Liebeck, G. M. Seitz,
\emph{Unipotent and Nilpotent Classes in Simple Algebraic
Groups and Lie Algebras},
Amer. Math. Soc. Providence, RI, (2012).

\bibitem[MaT]{MT} G. Malle and D. Testerman,
\emph{Linear Algebraic Groups and Finite Groups of Lie Type},
Cambridge Studies in Advanced Mathematics \textbf{133} (2011).



\bibitem[Sp2]{springer} T. A. Springer,
\emph{Linear Algebraic Groups. 2nd. edition},
Progress in Math. 9, Birkh\"auser, Boston (1998).

\bibitem[SS]{sp-st} T. A. Springer, R. Steinberg,
\emph{Conjugacy Classes, Seminar on Algebraic Groups and Related Finite Groups},
pp.167--266, Lect. Notes Math. 131, Springer,  (1970).

\bibitem[St1]{yale} R. Steinberg,
\emph{Lectures on Chevalley groups}, Yale University Press.



\bibitem[St2]{ihes} \bysame,
\emph{Regular elements of semisimple algebraic groups},
Inst. Hautes \'Etudes Sci. Publ. Math. No. 25 (1965) 49--80.

\bibitem[Su]{suzuki} M. Suzuki,
\emph{Group Theory I},
%Grundlehren der mathematischen Wissenschaften 247; 
Springer (1982).

\bibitem[TZ]{TZ}P. H. Tiep, A. E. Zalesskii,
\emph{Real conjugacy classes in algebraic groups and finite groups
of Lie type},
J. Group Theory \textbf{8} (2005) 291Ð315

\bibitem[W]{W} R. A. Wilson,
\emph{The Finite Simple Groups}. Graduate Texts in Math. 251. Springer.

\end{thebibliography}
\end{document}